\documentclass[reqno]{amsart}

\usepackage{latexsym}
\usepackage{amssymb}
\usepackage{palatino}
\usepackage{stmaryrd}
\def\x{\mbox{\boldmath $x$}}
\def\y{\mbox{\boldmath $y$}}
\def\z{\mbox{\boldmath $z$}}
\def\0{\mbox{\boldmath $0$}}

\def\cal{\mathcal}

\newtheorem{theorem}{Theorem}[section]
\newtheorem{lemma}[theorem]{Lemma}
\newtheorem{corollary}[theorem]{Corollary}

\newtheorem{definition}[theorem]{Definition}
\newtheorem{remark}[theorem]{Remark}

\begin{document}

%Version:D-G-\today

\title[Percolation in an ultrametric space]{Percolation in an ultrametric space}

\author[D.~A.~Dawson]{D.~A.~Dawson*}
\address{D.A. Dawson: School of Mathematics and Statistics, Ottawa, Canada K1S 5B6}
\email{ddawson@math.carleton.ca}
\thanks{* This research is  supported by NSERC}

\author[L.G. Gorostiza]{L.~G.~Gorostiza**}
\thanks{** Partially supported by CONACyT grant 98998}
\address{L.G. Gorostiza: CINVESTAV, Mexico City, Mexico}
\email{lgorosti@math.cinvestav.mx}



\today

\begin{abstract}
We study percolation in the hierarchical lattice of order $N$ where the probability of connection between two points separated by distance $k$ is of the form $c_k/N^{k(1+\delta)},\; \delta >-1$. Since the distance is an ultrametric, there are significant differences with percolation in the Euclidean lattice. We consider three regimes: $\delta <1$, where percolation occurs,  $\delta >1$, where it does not occur, and $\delta =1$ which is the critical case corresponding to the phase transition. In the critical case we use an approach in the spirit of the renormalization group method of statistical physics, and connectivity results of Erd\H{o}s-R\'enyi random graphs play a key role. We find sufficient conditions on $c_k$ such that percolation occurs, or that it does not occur. An intermediate situation called pre-percolation, which is necessary  for percolation, is also considered. In the cases of percolation we prove uniqueness of the constructed percolation clusters. In a previous paper \cite{DG1} we studied percolation in the $N\to\infty$ limit (mean field percolation), which provided a simplification that allowed finding a necessary and sufficient condition for percolation.  For fixed $N$ there are open questions, in particular regarding the behaviour at the   critical values of  parameters in the definition of $c_k$.
\end{abstract}

\maketitle

\noindent
\setcounter{section}{0}
\setcounter{equation}{0}
\noindent

\section{Introduction}
\label{S1}

Percolation theory in a lattice (e.g., the Euclidean lattice $\mathbb{Z}^d$) began with the work of Broadbent and Hammersley in 1957. The principal features of the model are that the space of sites is infinite and  its geometry plays an  essential role. The main problem is to determine if there is an infinite connected component, in which case it is said that percolation occurs. In the first models the connections (bonds) were  only between nearest neighbors.
(See \cite{G, Bo2} for background, and for a physics point of view, \cite{SA}.)
The question of percolation on noneuclidean graphs including nonamenable Cayley graphs given by finitely generated groups was formulated in \cite{BS}.
The study of long range percolation began in the mathematical physics literature (e.g., \cite{AN,  Be, NS, Sc}). In this case connections are allowed between points at any distance from each other with probability depending on the distance. The main problem is the same, and the geometry remains crucial. (See also \cite{BS,Bi,Bis,CGS,T}.)

The theory of random graphs started with the work of  Erd$^{{\rm ''}}\kern -.22cm$os and R\'enyi in 1959. The model consists of a finite number $n$ of vertices  with connection probability $p_n$ between pairs of vertices, depending on $n$ in some way, and there is no
structure on the set of vertices.  The results refer to what happens as $n\rightarrow \infty$, for example with the largest connected component. (See \cite{Bo, JLR, D} for  background.) Results from the theory of random graphs have been useful as technical tools in studies on percolation (e.g. \cite{BB, BLP}).

By introducing a structure on the set of vertices of a random graph or some special form of connection probabilities by means of a kernel which induces a sort of geometry, it is possible to generate a large class of interesting models which share properties of both percolation  and/or (classical) random graphs, and also ``small world'' random graphs (see e.g. the models in \cite{AB,Bo1,TV}).

On the other hand, hierarchical structures arise in the physical, biological and social sciences due to the multiscale organization of many natural objects (see e.g. \cite{BR,RTV}). In particular the hierarchical Ising model which was introduced by Dyson  \cite{Dy} has played an important role in statistical physics (see  \cite{BEI,CE,Si}) and in population genetics (see e.g. \cite{SF}).
Important applications of hierarchical structures have also been made by Kleinberg \cite{K, K2} in the area of search algorithms in computer science. A basic model is the {\it hierarchical group}  $\Omega_N$ of order $N$ (defined in Section 2), which can be represented as the set of leaves at the top of an infinite regular tree, where the distance between two points  is the number of levels from the top to their most recent common node. Such a distance satisfies the strong  triangle inequality
$$
d(x,y) \leq \max \{d(x,z), d(z,y)\}\quad \hbox{\rm for any}\quad x,y,z,
$$
which is the characteristic property of an ultrametric. (See e.g. \cite{Sch} for background on ultrametric spaces.) The main qualitative difference between  Euclidean-type lattices and an ultrametric space such as $\Omega_N$ is that in the former case it is possible to go far by a sequence of small steps, while in the latter case that is not possible, and the only way to go far is to make jumps of ever bigger sizes. This has important consequences for random walks for which there are analogies and differences with the Euclidean case (see e.g. \cite{DGW, DGW1} and references therein) and percolation in ultrametric spaces \cite{DG1}. In particular, percolation in $\Omega_N$ is possible only in the form of long range percolation, that is, with positive probabilities of connections between vertices separated by arbitrarily large distances.

With these precedents it is natural to investigate percolation in ultrametric spaces such as $\Omega_N$ where classical tools do not apply. Our aim is to develop  a mathematical framework that might be useful generally for this kind of model, not thinking about specific motivations from or applications to physics or any other field.

In \cite{DG1} we studied {\it asymptotic percolation} in $\Omega_N$ as $N\rightarrow \infty$
(or {\it mean field percolation}) with connection probabilities of the form $c_k/N^{2k-1}$ between two points separated by distance $ k$, and we obtained a necessary and sufficient condition for percolation. (See Subsection 3.1 for the definition of asymptotic percolation). The
Erd$^{{\rm ''}}\kern -.22cm$os-R\'enyi results on giant components of random graphs were a useful tool, although there are significant differences between classical random graphs and ultrametric ones (e.g., the average length of paths in the giant component of an ultrametric ball is much longer than in the classical case).

In the present paper we study percolation in $\Omega_N$ for fixed $N$ with connection probabilities of the form $c_k/N^{k(1+\delta)}$, $\delta >-1$, between two points separated by distance $k$. In this case percolation means that there is a positive probability that a given point of $\Omega_N$ belongs to an infinite connected component. This is a quite different situation from the asymptotic $N\rightarrow \infty$ model, and new methods must be used. However, properties of  Erd$^{{\rm ''}}\kern -.22cm$os-R\'enyi graphs are again useful, but now it is connectivity results that are of help,  specially the result in the Appendix based on Durrett's approach to connectivity \cite{D}. There are three regimes: $\delta <1$, $ \delta >1$ and $\delta =1$. Roughly speaking, under certain natural assumptions on $c_k$, for $\delta >1$ percolation does not occur, and for $\delta <1$ percolation occurs and the infinite connected component is unique. The most difficult is the critical case, $\delta =1$, where we use certain special forms of $c_k$ and percolation may or may not occur. Our aim is to find forms of $c_k$ such that percolation occurs, or that it does not occur, and in the case of percolation it turns out that the infinite connected component is unique.

We emphasize that taking the limit $N\to\infty$ provides a simplification which allowed us to obtain a sharp result, that is, necessary and sufficient condition for percolation in \cite{DG1}. In the case of finite $N$  that does not seem possible at present, and our work leads to some open problems, in particular regarding the behaviour at the  critical values of the parameters in the form of $c_k$ in the critical case $\delta=1$.  We also consider an intermediate situation  that we call {\it pre-percolation}, which is necessary  for percolation. Pre-percolation occurs in one of our results.

While we were working on this paper we learned about the manuscript of Koval et al \cite{KMT} (for which we thank them), where they also study percolation in $\Omega_N$ for fixed $N$, with connection probabilities of the form $1-\exp \{-\alpha /\beta^k\}$, $\alpha \geq 0, \beta>0$, between two points  separated by distance $k$. Some of their results for $\beta >1$ and  ours may be compared by setting $\beta =N^{1+\delta}$ (see Remark \ref{R3.2}).

For the cases where percolation occurs, specially with $\delta =1$, we use an approach in the spirit of the renormalization  group method of statistical physics  which has been employed, for example, in \cite{CE} for ferromagnetic systems on Dyson's hierarchical lattice $\Omega_2$   and for the study of long range percolation on $\mathbb{Z}^d$ (see  \cite{NS,Sc}).

In Section 2 we describe the model (the hierarchical group $\Omega_N$ and the associated random graph ${\cal G}_N$). In Section 3 we recall the result on mean field percolation  \cite{DG1} in order to compare it with a result in the present paper, and we state our   results for $\delta <1$ and $\delta >1$ (Theorem 3.1) and the critical case $\delta =1$ (Theorems 3.3 and 3.5), and we mention some open problems (Subsection 3.4).  Sections 4 and 5 contain the proofs. In an appendix we give a result  on connectivity of random graphs derived from \cite{D}, which is a key ingredient for the proof of percolation  in the critical case.

%\section{
\vglue .6cm
\noindent
\section{Description of the model}
%{\bf 2. Description of the model}
%\label{S2}
%\setcounter{section}{1}
\setcounter{equation}{0}
%\subsection
%\noindent{\bf 2.1 The hierarchical group $\Omega_N$}
\subsection{The hierarchical group $\Omega_N$}
\label{sub:2.1}

 For an integer $N\geq 2$, the {\it  hierarchical group of order $N$}, also called {\it hierarchical lattice of order $N$}, is defined as
\begin{eqnarray*}
\Omega_N &=&\{\mbox{\boldmath $x$} =(x_1, x_2,\ldots): x_i \in \{0,1,\ldots, N-1\}, i=1,2,\ldots, \sum_i x_i<\infty\}
\end{eqnarray*}
with addition componentwise mod $N$; in other words, $\Omega_N$ is a countable Abelian group given by the direct sum of a countable number of copies of the cyclic group of order $N$. The {\it hierarchical distance} on $\Omega_N$ is defined as
$$
d(\mbox{\boldmath $x$}, \y) =\left\{
\begin{array}{lcl}
0 & {\rm if} & \x=\y,\\
\max \{i: x_i \neq y_i\} & {\rm if} & \x \neq \y.
\end{array}\right.
$$
It is a translation-invariant metric which satisfies the strong (non-Archimedean) triangle inequality
$$
d(\x, \y)\leq \max \{d (\x, \z), d (\z, \y)\}\quad \hbox{\rm for any}\quad \x, \y, \z.
$$
Hence $(\Omega_{N}, d)$ is an ultrametric space, and it is well known that it can be represented as the leaves at the top of an infinite regular tree where $N$ branches emerge up from each node, and the distance between two points of $\Omega_N$ is the number of levels from the top to their most recent common node.

For each integer $k \geq 1$, a {\it $k$-ball} in $\Omega_N$, denoted  by $B_k$, is  a set of all points which are at distance at most $k$ from each other.
Any point of a ball can serve as a center. Once a center is chosen, one may speak of the interior and the boundary of the ball.
A $k$-ball contains $N^k$ points, and its boundary  contains $N^{k-1}(N-1)$ points. For $k>1$, a $k$-ball is the union of $N\;(k-1)$-balls, which are at distance $k$ from each other. For $j>k\geq 1$, a $j$-ball is the union of $N^{j-k}$ $k$-balls, which are at distance at least $k+1$ and at most $j$ from each other.
Two balls are either disjoint, or one is contained in the other (this is the reason why connections between nearest neighbours alone cannot produce percolation). For $k\leq j <\ell$, the ($j, \ell]$-{\it annulus} around $B_k$ is the set of all points $\y$ such that $j<d(\x, \y)\leq \ell$, where $\x$ is any point in $B_k$. The $(j, \ell]$-annulus is also described as $B_\ell \setminus B_j$ with $B_k \subset B_j \subset B_\ell$, and it contains $N^\ell(1-N^{j-\ell})$ points.  The number of points in a bounded subset $A\subset \Omega_N$ is denoted by $|A|$. The probability that a point chosen at random (uniformly) in $B_k$ belongs to $A\subset B_k$ is $|A|N^{-k}$.

We fix a point of $\Omega_N$ which we denote as $\bold{0}$.  Most of our considerations about percolation will refer to balls containing $\bold{0}$.

\subsection{The random graph ${\cal G}_N$}
\label{sub:2.2}
We define an infinite random graph ${\cal G}_N$ with the points of $\Omega_N$ as vertices, and for each $k\geq 1$ the probability of connection, $p_{\x,\y}$, between $\x$ and $\y$ with  $d(\x, \y)=k$ is given by
\begin{equation}
\label{eq:2.1}
p_{(k)}= p_{\x,\y}=\min\left(\frac{c_k}{N^{k(1+\delta)}},1\right),
\end{equation}
where $\delta >-1$ and $c_k >0$, all connections being independent. This can be realized in terms of a collection of independent uniform $[0,1]$ random variables
$\{ U_{(\x,\y)}\}$ by adding an (undirected) edge between $\x$ and $\y$ if and only if $U_{(\x,\y)}\leq p_{\x,\y}$  (see e.g. \cite{JLR}, page 4).

Our aim is to find sufficient conditions on $c_k$ and $\delta$ which imply that percolation occurs, or that it does not occur.

We will study separately  the  cases  $\delta >1$, $\delta <1$,  and   $\delta =1$. As we shall see, the  case $\delta = 1$ requires a more delicate analysis. In this case we take  $c_k$ of the following special forms:
\begin{itemize}
\item[(i)]
\begin{equation}
\label{eq:2.2}
c_k =C_0+C_1\log k+C_2k^\alpha,
\end{equation}
with constants $C_0\geq 0, C_1\geq 0, C_2\geq 0$ and $\alpha >0$.
\item[(ii)] We consider the  hierarchical distances with logarithmic scale
\begin{equation}
\label{eq:2.3}
k_n =k_n(K):=\lfloor K n\log n\rfloor, \quad n=1,2,\ldots,
\end{equation}
$K$ constant $ >0$, so that
\begin{equation}
\label{eq:2.4}
k_{n+1} -k_n \sim K\log n \quad {\rm as}\quad n\rightarrow \infty,
\end{equation}
(where $\sim$ has the usual meaning, see beginning of the proofs),
and we  consider the class of connection rates given by (\ref{eq:2.1}) with $\delta=1$, $c_k$ satisfying
\begin{eqnarray}
\label{eq:2.5}
c_{k_n} &=& C+a \log n \cdot N^{b \log n}
= C+a \log n \cdot n^{b \log N},
\end{eqnarray} and
\begin{equation}
\label{eq:2.6}
c_{k_n} \leq c_j \leq c_{k_{n+1}}\quad {\rm  for} \quad k_n<j <k_{n+1}.
\end{equation}
with constants $C\geq 0$, $a>0$, $b\geq 0$. The constant $K$ is chosen suitably in each case. The value $N=2$ is special because $\log 2  <1$,
and that is why for some results with $N=2$ we set $K> \frac{1}{\log 2}$.
\end{itemize}

The reason for considering $k_n(K)$ and  these forms for $c_{k_n}$ is that they provide  a four parameter family of comparison connection rates (with parameters $K,C,a,b$)  suitable for the renormalization analysis used in the proof
of percolation in Theorem \ref{T3.2}. An intuitive argument for this is given before Theorem \ref{T3.2}. Results of Theorem 3.5 are used to prove Theorem 3.3 which is our main result.

A  set of vertices any two of which are linked by a path of connections is called a {\it cluster}. By transitivity of $(\Omega_{N},d)$ we may focus on clusters containing ${\bf 0}$.
In the proofs  it is implicitly assumed that we consider a sequence of nested balls $(B_k)_{k\geq 1}$ such that
${\bf 0}$ belongs to $B_1$ (or to some $B_k$). We denote by $X_k$ the largest cluster contained in $B_k$, considering only connections within $B_k$ and not through points outside of $B_k$. If there are more than one largest cluster, we choose one of them uniformly from the existing ones. In this way each $k$-ball $B_k$ has a unique attached cluster $X_k$. Note that for $B_k \subset B_{k+j}$, either $X_k \cap X_{k+j}=\phi$ or $X_k \subset X_{k+j}$. Those clusters will be used only in the proofs of sufficient conditions for percolation.  The assumption of connections only within balls makes the renormalization approach quite practical for percolation, and connections through points outside would add to the possibility of percolation.

\begin{definition}
\label{D2.1}
{\rm We say that {\it percolation occurs} in ${\cal G}_N$ if there is a positive probability that a fixed point of ${\cal G}_N$ (for example
${\bf 0} $) belongs to an infinite cluster.}
\end{definition}

If percolation occurs,  the probability that there is an infinite cluster  is $1$. Indeed,  the event that there is an infinite cluster  is measurable with respect to the tail $\sigma$-algebra generated by the connections involving points outside each $k$-ball (containing $\0$) for every $k$, and the connections involving points  outside a ball are independent of  those inside, so by a $0\!$ - $\!1$ law   the probability that there is an infinite cluster is $0$ or $1$.

In some cases we consider percolation clusters of positive density, that is, $|X_k|N^{-k}$ does not decrease to $0$ as $k\to\infty$.

\begin{remark}\rm{
It follows immediately from the construction in terms of the family $\{U_{(\x,\y)}\}$ that given two families of connection probabilities $p^1_{\x,\y}$,
$p^2_{\x,\y}$ with $p^2_{\x,\y}\geq p^1_{\x,\y}$ for any $\x,\y$, percolation for family 1 implies percolation for family 2.
}
\end{remark}

%[Paragraph on the renormalization group approach].

%\section{\vglue .5cm
\noindent
%{\bf
\section{Results}
\label{S3}
\subsection{Mean field percolation}
\label{sub:3.1}
\setcounter{equation}{0}
For completeness, we start by recalling the result on  asymptotic percolation as $N\rightarrow \infty$ \cite{DG1}.  This will also be used for comparison with a result below.
The probability of connection between two points separated by distance $k$ is $c_k /N^{2k-1}$. Note that this corresponds to the critical case $\delta =1$ with $c_k$ in (\ref{eq:2.1}) multiplied by $N$, and this may be viewed as a normalization required for obtaining the  result in the limit. {\it Asymptotic percolation} is said to occur if
\begin{eqnarray*}
&& P_{{\rm perc}}: =\inf_{k\to\infty} \liminf_{N\rightarrow \infty} P({\bf 0}\;\hbox{\rm is linked by a path of  connections to a point at distance}\;k)>0.
\end{eqnarray*}
$P_{{\rm perc}}$ is the probability of percolation. For each $k\geq 1$, let $\beta_k\in (0,1)$ satisfy
\begin{equation}
\label{beta}
\beta_k =1-e^{-c_k \beta^2_{k-1} \beta_k}, \quad \beta_0=1,
\end{equation}
where $c_k \beta^2_{k-1}>1$. Note that $\beta_k$ is the well-known survival probability of a Poisson branching process with parameter $c_k\beta^2_{k-1}$. This corresponds to hierarchical level $k$, and the $\beta^2_{k-1}$ comes from the sizes of two connected giant components at the previous level $k-1$. Assume that $c_k \nnearrow \infty$ as $k\rightarrow \infty$, $c
_1>2\log 2$ and $c_2 >8 \log 2$. The results (see \cite{DG1}, Theorem 2.2 and Lemma 2.1) are that asymptotic percolation occurs if and only if
$$
\sum^\infty_{k=1} e^{-c_k}<\infty,
$$
and when it occurs, the probability of percolation is given by
$$
P_{\rm perc} =\prod^\infty_{k=1} \beta_k.
$$
(which is strictly positive if and only if the   exponential series converges), and percolation takes place through a cascade of clusters (in this case giant components) at consecutive hierarchical distances. For example, if $c_k=a \log k$ for large $k$, $a>0$, then asymptotic percolation occurs if and only if $a>1$. See Remark \ref{R3.3}(1) for a partially analogous result with fixed $N$.

\subsection{The  cases $\delta >1$ and $\delta <1$}
\label{sub:3.2}
\begin{theorem}
\label{T3.1}

\noindent
(a) If $\delta >1$ and $\sup_k c_k <\infty$, then  percolation does not occur.

\noindent
(b) If $\delta <1$ and $c=\inf_k c_k$ is large enough, then percolation occurs through a chain of clusters  in $k$-balls, and the percolation cluster is unique.
\end{theorem}

\begin{remark}\label{R3.2}
 {\rm Our results  and those of \cite{KMT} can be compared for $\beta >1$ therein, since in this case their connection probabilities $p_k=1-\exp(-\alpha/\beta^k)\sim\alpha/\beta^k$ as $k\to\infty$. If we set $\beta=N^{1+\delta}$ and let $c_k=c$ not depending on $k$ in (2.1), then  the decay rates agree  with $c$ corresponding  their $\alpha$. The comparison is between Theorem 1 in \cite{KMT} and our Theorem 3.1. We have that for $\delta>1,\beta>N^2$, percolation does not occur with any value of $c$, and for $-1<\delta<1,\; 1<\beta<N^2$, percolation occurs with $c$ sufficiently large, corresponding to $\alpha >\alpha_c(\beta)$ in \cite{KMT}. In addition in \cite{KMT} it is proved that there exists $\alpha_c(\beta)>0$ such that percolation does not occur for $\alpha <\alpha_c(\beta)$.
  Our main objective is to investigate the critical case $\delta=1,\beta=N^2$, for which  percolation does not occur for any $\alpha$ in \cite{KMT}. Our results in this case are stated in the next subsection. }
\end{remark}

\subsection{The case $\delta = 1$}
\label{sub:3.2x}

In the previous subsection we have seen that $\delta=1$ identifies the critical exponential decay rate for percolation. In this subsection we formulate our
main results  that determine the critical polynomial  rate for percolation.

\begin{theorem} \label{MR-1}    Let $\delta =1$ and
\[   c_k= C_0+C_1\log k+C_2k^\alpha \hbox{  with  }\alpha \geq 0,\]
where $C_0\geq 0$, $C_1\geq 0$, $C_2\geq 0$.

\noindent(a) If $\alpha >2$,
then for any $C_1$ there exist $C^*_0>0$ and $C_2^*>0$ such that if  $C_0>C_0^*$ and $C_2> C_2^*$,   percolation occurs and the percolation cluster is unique.

\noindent
(b) If $C_2=0$ and $C_1< N$, then percolation does not occur for any $C_0$.

\noindent
(c) If $\alpha >2$, there exists $C_*>0$ such that if $\max (C_0, C_1,C_2)<C_*$, then percolation does not occur.
\end{theorem}

The proof of this result is based on a renormalization argument that is formulated using the hierarchical distances $k_n(K)$ defined in (\ref{eq:2.3}) and the family of connection rates
$c_{k_n}$ defined in (\ref{eq:2.5}) with parameters $C,a,b$. This is the substance of Theorem 3.5.

In order to express one of the results below we introduce the following notion.

\begin{definition}\label{D3.1x}
{\rm We call {\em pre-percolation} the situation that (with probability $1$) there exists $n_0$ such that
    there is at least one connection from $(k_n,k_{n+1}]$ to $(k_{n+1},k_{n+2}]$ for all $n\geq n_0$.}
\end{definition}

Note that for percolation, in addition to pre-percolation, there would have to be paths connecting points in $(k_{n+1},k_{n+2}]$ which
are connected to $(k_n,k_{n+1}]$ to points in $(k_{n+1},k_{n+2}]$ which are connected to $(k_{n+2},k_{n+3}]$, etc.
\medskip

Before stating the next theorem, let us give an intuitive explanation for the choice of $k_n$ and $c_{k_n}$ in (2.3)-(2.6) with $b>0$, and the assumption $K<b$ (the assumption $\frac{2}{\log N}<K$ is a technical requirement for the method of proof). For this argument only we use the notation $\approx$ for approximate equality for large $n$ without giving it a rigorous meaning. The idea is that $k_n$ is the right scaling and the form of $c_{k_n}$ is the right
one which combines exactly with $k_n$ in order to produce the percolation cluster. We consider the largest clusters $X_{k_n}$ in each one of the $N^{k_{n+1}-k_n}\approx N^{K\log n}$ $k_n$-balls in a $k_{n+1}$-ball, and assume that their sizes are $|X_{k_n}|\approx \beta N^{k_n}$ for some $\beta\in (0,1)$, and that the probability of connection between two points in different clusters is $c_{k_n}/N^{k_{n+1}}$ (which is a lower bound for the actual probabilities). Let $s_n(\beta)$ denote the probability that two such clusters $X_{k_n}$ and $X^\prime_{k_n}$ in disjoint $k_n$-balls are connected. Then
\begin{eqnarray*}
&&s_n(\beta)\approx 1-\left(1-\frac{c_{k_{n}}}{N^{2k_{n+1}}}\right)^{|X_{k_n}||X^\prime_{k_n}|}\\
&&\approx 1-\left(1-\frac{c_{k_{n}}}{N^{2k_{n+1}}}\right)^{\beta^2N^{2k_n}}\\&&
\approx 1-\exp\left(1-\frac{c_{k_{n}\beta^2}}{N^{2(k_{n+1}-k_n)}}\right)\\
&&\approx \frac{C+a\log n\cdot N^{b\log n}}{N^{2K\log n}}\beta^2\\
&& \approx \frac{\beta^2a\log n}{N^{(2K-b)\log n}}=: r_n(\beta),\;{\rm with}\;\; b<2K,
\end{eqnarray*}
(see Lemma 5.8). Consider the E-R random graph
\[ G(N^{K\log n},r_n(\beta)),\]
and write $r_n(\beta)$ as
\[ r_n(\beta)= \frac{\beta^2a\log n}{N^{(K-b)\log n}}\frac{1}{N^{K\log n}}.\]
If $K>b$, then by the E-R theory only order $\log(N^{K\log n})$ of the $X_{k_n}$
are connected, hence the ratio of the size of the largest connected component in the $k_{n+1}$-ball to the size of the ball decreases to $0$
 as $n\to\infty$, so there cannot be a percolation cluster (of positive density).  Therefore we choose $K<b$. Now let $K_1=2K-b$, then $0<K_1<K$. Writing
 $r_n(\beta)$ as
 \[ r_n(\beta)=\frac{\beta^2 a}{K_1\log N}\log N^{K_1\log n}\frac{1}{N^{K_1\log n}}, \]
then by Theorem 2.8.1 of \cite{D} the probability that the graph $G(N^{K_1\log n},r_n(\beta))$ is connected tends to $1$ as $n\to\infty$ if $a>\frac{K_1\log N}{\beta^2}$. Connectivity of that graph whose vertices are the clusters $X_{k_n}$ in the $k_n$-balls in a $k_{n+1}$-ball means that the largest cluster in the $k_{n+1}$-ball contains the largest clusters in all the $k_n$-balls it contains.  If this can be proved for all sufficiently large $n$, then percolation follows.
Note that this argument provides a percolation cluster if $K<b$, but it does not imply that percolation does not occur if $b>K$.

 \medskip
\begin{theorem}
\label{T3.2} Let $\delta =1$.

\noindent
%(a) Let $c_k =Ck^\alpha$ (see (2.2)). If $\alpha <-1$, then  percolation does not occur.

\noindent
(a) Let  $C_2=0$ in (\ref{eq:2.2})  and  $c_{k_n}=C+aN\log n$ in (\ref{eq:2.5}), with $K=1$ and $N\geq 3$  in (\ref{eq:2.3}).
%(but this value of $K$ is not essential).
Then

(i)  for $a<1$,  percolation does not occur,

(ii)  for $a >1$, pre-percolation occurs.
\medskip

\noindent
(b) Assume that $\{c_{k_n}\}$ satisfy (\ref{eq:2.5}), (\ref{eq:2.6}) with $b>0$ and  $c_{k_n}=C+a\log n\cdot N^{b \log n}$ where $k_n=k_n(K)$ is given by (\ref{eq:2.3}) and the pair $(K,b)$ satisfy
\begin{equation*}
\frac{2}{\log N}<K<b.
\end{equation*}

%(iii) for $N=2$, let $b>\frac{3}{\log 2}$, and $K>\frac{1}{2}\left(\frac{1}{\log 2}+b\right)$.

(1) Then there exist $C>0$ and $a_*>0$ such that for $a>a_*$ there is a sequence $(\beta_{n})_n$ such that
\begin{equation}
\label{eq:3.1}
\liminf_{n\rightarrow \infty} \beta_{n}>0,
\end{equation}
and for the clusters $X_{k_n}$ in a nested sequence of $k_n$-balls $B_{k_n}$ containing  ${\bf 0}$,
%\begin{equation}
%\label{eq:3.2}
%\liminf_{n\rightarrow \infty}  P(|X_{k_n}|\geq \beta_{n}N^{k_n})>0,
%\end{equation}

\begin{equation}
\label{eq:3.2}
  P(\hbox{there exists }\; n_{00}\hbox{  such that  }|X_{k_n}|\geq \beta_{n}N^{k_n}\hbox{  for all }n\geq n_{00})=1,
\end{equation}
 percolation occurs, and the percolation cluster is unique.
\medskip

(2)  Assume that $b<  2K-\frac{1}{\log N}$. Then the percolation cluster is given by a ``cascade'' of clusters at distances $k_n$, more precisely,  there exists a (random) number $n_0$ such that for $n\geq n_0$ connections between $X_{k_n}\cap (B_{k_n}\backslash B_{k_{n-1}})$ and $X_{k_{n+\ell}}\cap (B_{k_{n+\ell}}\backslash B_{k_{n+\ell-1}})$ occur only for $\ell=1,2$.

\noindent
(c) In the special case with connection probabilities $c_j=\frac{c_{k_{n}}}{N^{2k_{n+1}}}$ for $k_n+1\leq j\leq k_{n+1}$ and $0<b\leq \frac{2}{\log N} <K$, percolation does not occur.
\end{theorem}

\medskip

\begin{remark}
\label{R3.3}

\noindent
(1) {\rm In Theorem \ref{T3.2} (a) a slightly different result holds for $N=2$ since $\log 2 <1$ so that we need $K>\frac{1}{\log 2}$.}

\noindent(2) {\rm  Note the consistency of Theorem \ref{T3.2} (a) with the example of asymptotic percolation recalled in Subsection \ref{sub:3.1}.  The difference is that  in the finite $N$ case we only have pre-percolation, and  the connections are at hierarchical distances $k \log k$ rather than $k$.}

\noindent
(3) {\rm Theorem \ref{T3.2}(a)(ii) implies that in part (b) pre-percolation occurs with any $a>0$ and $b>0$.}

\noindent
%(3) {\rm It turns out that pre-percolation is necessary  for percolation.  This follows from Lemma 5.1(f) for $b=0$ and from Lemma 5.2 for $b>0$.}

\noindent
%(4) {\rm To see that pre-percolation does not imply percolation, let $\{c_{k_n}\}$ be given by (\ref{eq:2.5}) with $a>0$,  $0<b<1$ and $K=1$. Then %pre-percolation
%occurs (Remark \ref{R3.3} (2)). By an argument similar to that before Theorem 3.5 (but now rigorous), since $|X_{k_n}|\leq N^{k_n}$ then the probability  $s_n$ %of connection between two clusters in disjoint $k_n$-balls in a $k_{n+1}$-ball satisfies
%\[ s_n\leq \frac{C+a\log n\cdot N^{b\log n}}{N^{\log n}}\frac{1}{N^{\log n}}\]
%for all sufficiently large $n$. Since
%\[ \frac{C+a\log n\cdot N^{b\log n}}{N^{\log n}}\to 0 \;{\rm as  }\; n\to \infty,\]
%then by the E-R theory there is no percolation for this choice of the parameters $(a,b,K)$.
% (Note however that this argument is for a fixed scale $k_n$ determined by the parameter $K$. For any $0<b<1$ and $N$ such
%that $b>\frac{2}{\log N}$, if $K$ is chosen so as to satisfy the hypothesis of Theorem 3.3(c), then percolation can be achieved.)
%}

\noindent
(4) {\rm The cascade of   clusters in Theorem \ref{T3.2} (b)(2) is analogous to the  cascade of giant components in the  mean field case  \cite{DG1}.}

\noindent
(5) {\rm The formulation in Theorem \ref{T3.2}(b) is used as a technical tool to prove the  result in Theorem \ref{MR-1}, and also provides  a setting to
give a refined result.}

\end{remark}

\subsection{Open problems and related developments.}
\label{R3.4}

\noindent{\em (1)}
In Theorem 3.5(a) we have proved that with $c_{k}=C_{0}+C_{1}\log
k+C_{2}k^{\alpha}$ percolation does not occur if $C_{2}=0$ and $C_{1}%
<{N}$. On the other hand in (b) we proved that for $\alpha>2$ and
$C_{0}$ sufficiently large$,$ percolation occurs. It remains an open question
as to whether percolation can occur for all $\alpha>0$ or even for $C_2=0$
and some $C_{1}$ sufficiently large$.\;$We next explain that to resolve these
questions analogues of well-known results for Erd\H{o}s-R\'enyi graphs would be needed for
a class of  ultrametric random graphs. An  {\em ultrametric random graph} URG(M,d) is a random graph on a finite
set $M$ with ultrametric $d$ and with  connection probabilities $p_{x,y}$ that depend on the ultrametric distance $d(x,y)$.

\medskip

Consider the case $C_{2}=0$. The expected number of in-edges
to the annulus $(k_{n+1}, k _{n+2}]$ from $B_{k_{n+1}}$ is of order $
O(a(1-\frac{1}{N})\log n)$
and the expected number of out-edges to the $(k_{n+2}, k _{n+3}]$ annulus is of order\newline
$O(a(1-\frac{1}{N})\log(n+1)).$ \ In order to determine  the number of in-edges that connect to
an out-edge, it would be necessary to determine the probability that two
randomly chosen vertices in the $(k_{n+1},k_{n+2}]$-annulus are connected by a path in the
associated ultrametric random graph. \ This is
related to the problem of determining the distribution of sizes of the connected components. \ These are
 open problems.

\medskip

Consider the case $C_{2}>0$. In Theorem 3.5(c)  we have proved that the random
graph based on a lower bound for the connection probabilities at distances
$k_{n}+1,\dots,k_{n+1}$ does not exhibit percolation in the case
$0<b\leq\frac{2}{\log N}<K$ \ (corresponding to the case $\alpha<2$). \ In order
to refine the argument and determine the behaviour for the actual connection
probabilities which arise if $0<\alpha<2$ it would be necessary to
determine the size of the largest connected component, that is, the number of
$k_{n}$-balls (more precisely their largest connected components) in the
$k_{n+1}$-ball which are connected and in the largest connected component in
the associated ultrametric random graph  (as $n\rightarrow\infty)$. This is an open problem.

% In Theorem \ref{T3.2} we considered the cases
%\begin{equation} \label{OP1} c_{k_n}=C+a\log n
%\end{equation}
%and

%\begin{equation} \label{OP2}  c_{k_n}=C+a\log n \cdot N^{b\log n}.
%\end{equation}

%In case (\ref{OP1}) we have established non-percolation for $a<1$ and pre-percolation for $a>1$.
%It remains an open question if  percolation can occur for $a>$ some $a^*>1$.

%Part (c) implies that in case (\ref{OP1}) percolation does not occur for any $a$.
%In case (\ref{OP2}), the question as to whether given parameters $a,b,K$ satisfying the conditions of Theorem \ref{T3.2} non-percolation can occur for some
%$C>0$ or  $0<a<a_*$ is open.
 %and, since pre-percolation occurs with any $a>0$ and $b>0$ (Remark 3.5(2)), the question of an additional condition for percolation is also open; a partial
 %answer is given by Theorem 3.3(c).

% \textbf{REMOVE OR REVISE} If those questions for either case have positive answers, there would be corresponding critical values of $a$ (presumably $a=1$ in %case (\ref{OP1})), and the behaviour at the critical value would be of interest. { The value of $a_*$ in Theorem \ref{T3.2}(c) depends on  lower bounds that are %conveniently chosen in the proof.  In principle they can be changed to obtain lower values of $a_*$ but it remains to find the optimal value of $a_*$.} %Increasing the values of $a$ and$/$or $b$ favours percolation and decreasing them inhibits percolation. The question of a critical value for $b$ for percolation %is also open.

\medskip
\noindent
{\em(2)} It would also be of interest to consider the intermediate case with $\delta=1$ but with connection probabilities of the form
 $\frac{n^{(\log n)^k}}{N^{2k_n}}$.

\medskip

 \noindent
{\em(3)}  One can ask if there can exist a zero-density infinite cluster for parameter values where positive density connected clusters do not exist.

%\textbf{XXX: Do we still want this?} The unique percolation clusters obtained in Theorems 3.1(b) and 3.3(c) are dense in the sense that the contributing clusters %$X_{k_n}$ occupy a sizeable proportion of the balls $B_{k_n}$ they belong to. One can ask if there could exist other ``thin'' (zero-density) percolation clusters %not intersecting the dense ones or if these can occur at the critical point (if it exists).  The existence of a non-uniqueness phase of nearest neighbour %percolation on nonamenable Cayley graphs was conjectured by Benjamini and Schramm \cite{BS} and established by Pak and Smirnova-Nagnibeda \cite{PSN}. It would be %of interest to determine if this can occur in our case.
\medskip

\noindent
%(3)  In Theorem \ref{T3.2}(a) we don't know what happens with $\alpha = -1$.

\noindent
{\em (4)} Berger \cite{Be} has studied the behaviour of  random walk on the infinite cluster of long-range percolation in Euclidean lattices of dimensions $d=1,2$.  It would also be interesting to investigate this behaviour on the infinite clusters obtained in Theorems 3.1(b) and  \ref{T3.2} (b). Long-range random walks on $\Omega_N$ have been studied in \cite{DGW}.

\subsection{The renormalization group approach}

%\textbf{Shorten further ??}

%\textbf{Shorten or delete this subsection ??}

%\bigskip

The basic strategy we employ is in the spirit of the renormalization group method of statistical  physics \cite{CE}, which has been used by Newman and Schulman \cite{NS}, Section 2,  in their study of long range percolation in the Euclidean lattice.

Consider the countable ultrametric space $(\Omega_N,d)$.  For each integer $k $ we define an equivalence relation on $\Omega_N$ by
$$ \x\equiv_k \y \hbox{ iff }d(\x,\y)\leq k,\; \hbox{that is, } \x \hbox{ and }\y\hbox { belong to the same }k\hbox{-ball}.$$
Now consider the set of equivalence classes $\mathcal{E}_k$ furnished with the ultrametric

$$d_k(\widetilde {\x}, \widetilde{\y}):= d(\x,\y)=0\hbox{ if } d(\x,\y)\leq k,\\
d_k(\widetilde {\x}, \widetilde {\y}):= d(\x,\y)-k\hbox{ if } d(\x,\y)> k,$$
where $\widetilde {\x}$ is an equivalence class containing a point $\x$.
Then the resulting set of equivalence classes with ultrametric $d_k$ can be identified with
$(\Omega_N,d)$.

Given a graph
$$
{\cal G}_N= (\Omega_N, {\cal E}_N),
$$
with edges $\cal E_N$   given by  a symmetric subset of $\Omega _N\times \Omega_N$
 we obtain a new graph  as follows. The set of vertices is the set of all $k$-balls, and the set of edges $\mathcal{E}_{N,1}$
are such that $\widetilde {\x}$ has a connection to $\widetilde{ \y}$ if there is a connection in $\mathcal{G}_N$ between the largest connected subset (cluster) of $\widetilde {\x}$ and the largest connected subset of $\widetilde {\y}$. Using the above identification this defines a new graph $\mathcal{G}^1_N=\Phi\mathcal{G}_N$  on $(\Omega_N,d)$.  Iterating this procedure we obtain a sequence
$${\cal G}^k_N= \Phi^k\mathcal{G}_N,\;\; k\geq 1,$$ of graphs all having vertex set $\Omega_N$.  In addition we assign to each vertex $v$ in ${\cal G}^k_N$ the $[0,1]$-valued random variable  $$Y_{k}(v):=\frac{|X_{k}(v)|}{N^{k}}$$ where $X_{k}(v)$ denotes the set of vertices in  $\mathcal{G}_N$  contained in the  connected cluster  in the $k$-ball corresponding to $v$ obtained as the union of the clusters in the $(k-1)$-balls it contains.
This construction  defines the renormalization mapping $\Phi: {\cal G}_N\rightarrow {\cal G}_N$ such that $\Phi^k :{\cal G}_N \rightarrow {\cal G}^{k}_N$ for each $k$. Note that at each iteration the connection probabilities are the probabilities that the largest connected subsets of equivalence classes are connected and that these probabilities are random and dependent (because the connected components have random sizes), and they change at each iteration.

Rather than working directly with the sequence ${\cal G}^{k}_N$ we choose a subsequence $k_n$ and construct a sequence of renormalization maps $\Phi^{k_n}$ such that the number of points in a ball of radius 1 (with respect to the new distances $d_{k_n}$) increases to infinity as $n\to\infty$.
In particular, we will show that there exists an increasing  sequence of integers $(k_n)$ (see (2.3)) and a sequence of graphs $\mathcal{G}^{k_n}$ with $N^{k_n-k_{n-1}}$ vertices constructed recursively as follows, that is,

$$\mathcal{G}^{k_{n+1}}=\Phi_n \mathcal{G}^{k_{n}},$$
where $\Phi_n$ depends on $n$ since it is a mapping from a graph with vertices $\Omega_{N^{k_n-k_{n-1}}}$ to a graph with vertices $\Omega_{N^{k_{n+1}-k_n}}$ and with connection probabilities between vertices that are a function of the distance between them.
 Moreover we can identify $\mathcal{G}^{k_n}$ with a subgraph of $\mathcal{G}_N$ and these subgraphs are a decreasing function of $n$.  We establish percolation by showing that the intersection of these subgraphs starting at a given  point in $\Omega_N$ is non-empty with positive probability.
 The difference now is that $\Phi^{k_{n+1}}$ is not obtained by iteration but by means of $\Phi_n$:
$$ \Phi^{k_{n+1}}=\Phi_n\Phi^{k_n}.$$

Given the sequence $(k_n)$ we can consider the equivalence classes given by

$$ \x\equiv_{k_n} \y\;\; \hbox{iff} \;  d(\x,\y) \leq k_n,$$
and define the ultrametric  $d_{k_n}$ by
%$$ d_{k_n}(\widetilde {\x},\widetilde {\y})= 1\hbox{  iff  } k_n<d(\x,\y)\leq k_{n+1},
%$$
%$$ d_{k_n}(\widetilde {\x},\widetilde {\y})= 2\hbox{  iff  } k_{n+1}<d(\x,\y)\leq k_{n+2},
%$$
$$ d_{k_n}(\widetilde {\x},\widetilde {\y})= \ell\hbox{  iff  } k_{n+\ell-1}<d(\x,\y)\leq k_{n+\ell}, \; \ell\geq 1.$$
We now consider the graph $\mathcal{G}^{k_n}$ whose vertices are the  $d_{k_n}$ equivalence classes.
Two points in $\mathcal{G}^{k_n}$  at $d_{k_n}$-distance $\ell \geq 1$ are connected if there is a $\mathcal{G}_N$-edge joining the largest connected components in these equivalence classes in $\mathcal{G}^{k_{n-1}}$.
Note that there  are $N^{k_{n+1}-k_n}$ points in a ball of $d_{k_n}$-radius 1, $N^{k_{n+2}-k_n}$ in a ball of radius 2, etc. The proof of percolation in the case $\delta <1$ given in Section 4 involves showing that as $n\to\infty$ the graphs $\mathcal{G}^{k_n}$ occupy a certain portion of the $k_n$-balls and an increasing sequence can be linked in a cascade with probability approaching $1$.

In order to apply these ideas to the more delicate critical case $\delta =1$  in Section 5 we define $Y_{k_n}(v),\;v\in \mathcal{G}_{k_n}$ as above.
Then given the random graph $\mathcal{G}_N$, the nontriviality, $\liminf_{n\to\infty}Y_{k_n} >0$,  has probability 0 or 1.  Our goal is to find a sufficient condition for this to be 1.
%We again consider the mappings
%$$ F^N_{k_n}\to F^N_{k_{n+1}}=\widehat \Phi_{k_n}F^N_{k_n}.$$
%Since we no longer take the limit $N\to\infty$, the mappings do not reduce to mappings on $[0,1]$ as in the mean field case.  However since  the random graph %$\mathcal{G}^{k_n}$ has $N^{{k_n}-k_{n-1}}$ vertices and this goes to infinity as $n\to\infty$, we expect the asymptotic behaviour to be close to deterministic.
In order to achieve this our strategy  is to look for a pair of sequences
$k_n\to\infty$,  $\liminf_{n} \beta_n >0$, such that the probability that $Y_{k_n}\geq \beta_n$ converges to $1$ as $n\to\infty$.  This program will be carried out in Subsection 5.2 using as basic tools large deviation estimates for binomial distributions and probability bounds for the connectivity of an Erd\H{o}s-R\'enyi graph.

\section{Proofs for the  cases $\delta>1$ and $\delta <1$}
\label{S4}
\setcounter{equation}{0}

We first mention a few notational points.  In some  places in the proofs in this and the following section where a number appears which should be a non-negative integer and it is not necessarily so, it should be interpreted as its integer part.

 $a_n\sim b_n$ means that $\frac{a_n}{b_n}\to 1$ as $n\to\infty$, $a_n>>b_n$ means that $b_n=o(a_n)$,  and \[a_n \lesssim b_n {\rm \;means\; that\;}
0\leq \liminf_n \frac{a_n}{b_n}\leq \limsup_n \frac{a_n}{b_n} \leq 1.\]

\begin{definition}
\label{D4.1}  {\rm
For $0<\gamma <1$, we say that a $k$-ball} $B_k$ is $\gamma$-good
{\rm if its attached cluster} $X_k$  {\rm satisfies} $|X_k|\geq N^{\gamma k}$.
\end{definition}

If ${\bf 0} \in B_k$ and  $B_k$ is $\gamma$-good,  then the probability that ${\bf 0} \in X_k$ is greater than or equal to  $ N^{(\gamma -1)k}$.

\bigskip

\noindent
{\bf Proof of Theorem 3.1.}

%\begin{itemize}
%\item[\kern -1cm (a)]
\noindent
(a) If  suffices to show that
\begin{equation}
\label{eq:4.1'}
P(B_k\;\hbox{\rm is connected to its complement for infinitely many}\;k)=0.
\end{equation}

For $j\geq k$,
\begin{eqnarray*}
\lefteqn{P(B_k\;\hbox{\rm is connected to the}\;(j,j+1]
\hbox{\rm -annulus around it})}\\
&=& 1-\biggl(1-\frac{c_{j+1}}{N^{(j+1)(1+\delta )}}\biggr)^{N^k N^j (N-1)}\\
&&(\hbox{\rm using}\; 1-(1-x/y)^z<xz/y, \quad 0<x<y, \quad z\geq 2)\\
&<& c_{j+1}\frac{N^kN^j (N-1)}{N^{(j+1)(1+\delta )}}\\
&\leq & M\frac{N^k}{N^{j\delta}},\quad M\;\hbox{\rm is a constant }(\hbox{since }\sup_k c_k<\infty),
\end{eqnarray*}
hence
$$
P(B_k \;\hbox{\rm is connected to its complement)}\;\leq M
\sum_{j\geq k}\frac{N^k}{N^{j\delta}},
$$
and since
$$
\sum_{k\geq 1} \sum_{j\geq k} \frac{N^k}{N^{j\delta}}<\infty,
$$
(\ref{eq:4.1'})  follows by Borel-Cantelli.

\noindent
(b) For each $n$, let $k_n$ be given by:
\begin{equation}
\label{eq:4.1}
k_n =\lfloor n\log n\rfloor ,
\end{equation}
hence
\begin{equation}
\label{eq:4.2'}
k_{n+1} -k_n \sim \log n \quad{\rm as}\; n\rightarrow \infty.
\end{equation}
Choose $\gamma$ so that
\begin{equation}
\label{eq:4.2}
\frac{1+\delta}{2}<\gamma <1.
\end{equation}
Then the probability of connection between two  clusters $X_{k_n}$ and $X'_{k_n}$ in (disjoint) $\gamma$-good $k_n$-balls in  a $k_{n+1}$-ball  is bounded below by
\begin{equation}
\label{eq:4.3'}
 1-\biggl(1-\frac{c}{N^{k_{n+1} (1+\delta)}}\biggr)^{|X_{k_n}||X'_{k_n}|} \geq 1-\biggl(1-\frac{c}{N^{k_{n+1}(1+\delta)}}\biggr)^{N^{2\gamma k_n}}.
\end{equation}

There are $N^{k_{n+1}-k_n}\;\; k_n$-balls in a $k_{n+1}$-ball. If at least $N^{\gamma (k_{n+1}-k_n)}$ of  them are $\gamma$-good and their  clusters are connected within the $k_{n+1}$-ball, then the size of the  cluster   $X_{k_{n+1}}$ in the $k_{n+1}$-ball is
 greater than or equal to $ N^{\gamma k_{n+1}}
$, so the $k_{n+1}$-ball is $\gamma$-good.

For each $n$, let $p_n$ denote the probability that $B_{k_n}$ is $\gamma$-good:
\begin{equation}
\label{eq:4.3}
p_n =P(|X_{k_n}|\geq N^{\gamma k_n}).
\end{equation}
Then, writing ${\rm Bin} (n, p, \geq k)=\sum^n_{j=k}({n\atop j})p^j (1-p)^{n-j}$ and ${\rm Bin}(n, p, <k)=1-{\rm Bin}(n,p, \geq k)$, we have, by (\ref{eq:4.3'}), (\ref{eq:4.3}), and independence of connections at different hierarchical distances,
\begin{eqnarray}
\label{eq:4.4}
&& p_{n+1} \geq  {\rm Bin} \biggl(N^{k_{n+1}-k_n},p_n, \geq N^{\gamma (k_{n+1}-k_n)}\biggr)
 \biggl[1-\biggl(1-\frac{c}{N^{k_{n+1}(1+\delta)}}\biggr)^{N^{2\gamma k_n}}\biggr]^{N^{\gamma (k_{n+1}-k_n)}}.
\end{eqnarray}
Now,
\begin{eqnarray}
\label{eq:4.4''}
\lefteqn{1-\biggl(1-\frac{c}{N^{k_{n+1}(1+\delta)}}\biggr)^{N^{2\gamma k_n}}}\nonumber\\
&\geq & 1-\biggl(1-\frac{cN^{2\gamma k_n -(1+\delta)k_{n+1}}}{N^{2\gamma k_n}}\biggr)^{N^{2\gamma k_n}}\nonumber\\
&&(\hbox{\rm using}\; (1-x/y)^y <e^{-x},\quad x>0,\, y>0),\nonumber\\
&>& 1-\exp \{-c N^{2\gamma k_n-(1+\delta)k_{n+1}}\},
\end{eqnarray}
\noindent
hence
\begin{eqnarray}
\label{eq:4.4'}
&&\biggl[1-\biggl(1-\frac{c}{N^{k_{n+1}(1+\delta)}}\biggr)^{N^{2\gamma k_n}}\biggr]^{N^{\gamma (k_{n+1}-k_n)}}\nonumber\\
&&>(1-\exp \{-c N^{2\gamma k_n-(1+\delta)k_{n+1}}\})^{N^{\gamma (k_{n+1}-k_n)}}.
\end{eqnarray}
(\ref{eq:4.1}) and (\ref{eq:4.2}) imply that
\begin{equation}
\label{eq:4.5}
2\gamma k_n-(1+\delta) k_{n+1} > \varepsilon n\log n
\end{equation}
for some $0<\varepsilon <1$ and  sufficiently large $n$. Then, from
(\ref{eq:4.4'}),
\begin{equation}
\label{eq:4.6}
\biggl[1-\biggl(1-\frac{c}{N^{k_{n+1}(1+\delta)}}\biggr)^{N^{2 \gamma k_n}}\biggr]^{N^{\gamma (k_{n+1}-k_n)}} > 1-\varepsilon_n,
\end{equation}
where
\begin{equation}
\label{eq:4.7}
\varepsilon_n :=1-\biggl(1-\exp \{ -c N^{2\gamma k_n -(1+\delta )k_{n+1}}\}\biggr)
^{N^{ \gamma (k_{n+1}-k_n)}},
\end{equation}
and we have that
\begin{equation}
\label{eq:4.8}
\sum_n \varepsilon_n <\infty ,
\end{equation}
because, from (\ref{eq:4.7}), (\ref{eq:4.2'}),   and (\ref{eq:4.5}),
\begin{eqnarray*}
\lefteqn{
\varepsilon_n < N^{\gamma (k_{n+1}-k_n)}\exp \{-c N^{2\gamma k_n -(1+\delta ) k_{n+1}}\}}\\
&\lesssim & N^{\gamma \log n}\exp \{-cN^{\varepsilon n\log n}\}\\
&=&n^{\gamma \log N}\exp \{-cn^{\varepsilon n\log N}\}
\end{eqnarray*}
for all sufficiently large $n$, which is summable.

From (\ref{eq:4.4}) and (\ref{eq:4.6}),
\begin{equation}
\label{eq:4.9}
p_{n+1}> {\rm Bin} (N^{k_{n+1}-k_n}, p_n, \geq N^{\gamma (k_{n+1}-k_n)}) (1-\varepsilon_n).
\end{equation}

Now we assume that there exists $n_0$ such that
\begin{equation}
\label{eq:4.10}
 N^{k_{n_0+1}-k_{n_0}} p_{n_0}>>N^{\gamma (k_{n_0+1}-k_{n_0})},
\end{equation}
(this assumption will be verified below), and we suppose that $p_{n_0}>\eta >0$. Then,  by a large deviation inequality \cite{JLR} (p. 26, (2.6)),
%\cite{JLR[p.26, (2.6)},
\begin{eqnarray}
\label{eq:4.11}
&& {\rm Bin} \left(N^{k_{n_0+1}-k_{n_0}}, p_{n_0}, \geq N^{\gamma (k_{n_0+1}-k_{n_0})}\right)\\
&&= 1-{\rm Bin} \left(N^{k_{n_0+1}-k_{n_0}}, p_{n_0}, <N^{k_{n_0+1}-k_{n_0}}p_{n_0}
-(N^{k_{n_0+1}-k_{n_0}}p_{n_0} -N^{\gamma (k_{n_0+1} -k_{n_0})})\right)\nonumber\\
&&\geq  1-\exp
\biggl\{-\frac{(N^{k_{n_0+1}-k_{n_0}} p_{n_0} -N^{\gamma (k_{n_0+1}-k_{n_0})})^2}{2N^{k_{n_0+1}-k_{n_0}}p_{n_0}}\biggr\},\nonumber
\end{eqnarray}
hence, by (4.15),
\begin{equation}
\label{eq:4.12}
{\rm Bin} \left(N^{k_{n_0+1}-k_{n_0}}, p_{n_0}, \geq N^{\gamma (k_{n_0+1}-k_{n_{0})}}\right)\geq 1-\delta_{n_0},
\end{equation}
where we may
 take
\begin{equation}
\label{eq:4.13}
0<\delta_{n_0} <(1-\varepsilon_{n_0}-\eta) /(1-\varepsilon_{n_0})
\end{equation}
(with $\varepsilon_{n_0}$ and $\eta$ small enough so that $\varepsilon _{n_0}+\eta <1)$. Then, from (\ref{eq:4.9}), (\ref{eq:4.12}) and (\ref{eq:4.13}),
\begin{equation}
\label{eq:4.14}
p_{n_0+1}> (1-\delta_{n_0} )(1-\varepsilon_{n_0})>\eta,
\end{equation}
so, iterating the argument,
\begin{equation}
\label{eq:4.14'}
p_n >\eta \quad \hbox{\rm for all}\quad n\geq n_0,
\end{equation}
and therefore (\ref{eq:4.10}) will in fact hold for all sufficiently large $n_0$.

From (\ref{eq:4.14}),
\begin{equation}
\label{eq:4.15}
1-p_{n_0+1}\leq \delta _{n_0} +\varepsilon_{n_0}.
\end{equation}
Choosing $\delta_{n}, n\geq n_0$, so that $\sum\limits_{n\geq n_0} \delta_{n}<\infty$, we then have from (\ref{eq:4.15}), together with (\ref{eq:4.8}),
$$
\sum_{n\geq n_0} (1-p_n)<\infty.
$$
Therefore, by Borel-Cantelli,
\begin{equation}
\label{eq:4.22}
P(\hbox{\rm all but finitely many}\;B_{k_n}\;{\rm are}\;\gamma\hbox{\rm -good})=1.
\end{equation}

On the other hand, the probability that the clusters $X_{k_n}$ and $X_{k_{n+1}}$ of  the $\gamma$-good  balls $B_{k_n}$ and
$B_{k_{n+1}}, B_{k_n} \subset B_{k_{n+1}}$, are not  connected is, for large $n$, by (\ref{eq:4.1}) and  (\ref{eq:4.5}), bounded above by
\begin{eqnarray*}
\biggl(1-\frac{c_{k_{n}}}{N^{k_{n+1} (1+\delta)}}\biggr)^{|X_{k_n}||X_{k_{n+1}}|} &\leq & \biggl( 1-\frac{c}{N^{k_{n+1}(1+\delta)}}\biggr)^{N^{\gamma (k_n +k_{n+1})}}\\
&<& \exp \{-c N^{\gamma (k_n+k_{n+1})-k_{n+1}(1+\delta)}\}\\
&\lesssim & \exp \{-c N^{(2\gamma k_n-(1+\delta )k_{n+1})}\} \\
&<&\exp \{-cN^{\varepsilon n \log n}\}\\
&=& \exp \{-cn^{\varepsilon n\log N}\},
\end{eqnarray*}
which is summable. Therefore, by Borel-Cantelli,
\begin{equation}
\label{eq:4.23}
P \hbox{\rm (all but finitely many pairs}\; X_{k_n}\;{\rm and}
\; X_{k_{n+1}}\;\hbox{\rm  are connected}) =1.
\end{equation}

Since there is a positive probability that $\0$ belongs to some $X_{k_n}$ (the probability is greater than or equal to $ N^{(\gamma-1)k_n})$ in the  chain of $\gamma$-good  $B_{k_n}$'s with connected clusters $X_{k_n}$, then by (\ref{eq:4.22}) and (\ref{eq:4.23})  there is an infinite cluster  that contains $\0$ with positive probability. So, percolation occurs.

It remains to verify assumption (\ref{eq:4.10}). We do this by means of a connectivity result for E-R random graphs. Writing $c_k/N^{k(1+\delta)}$ in (\ref{eq:2.1}) as $\widetilde{c}_k/N^k, \widetilde{c}_k = c_k /N^{\delta k}$, we have that for all $k\leq n$,
$$
\widetilde{c} _k > \overline{c}_n \log N^n\quad {\rm with}\quad \overline{c}_n =\frac{c}
{N^{ \delta n} \log N^{n}} .
$$
We consider the E-R random graph $G(N^n, \overline{c}_n/N^n)$ whose vertices are the points of an $n$-ball. If
\begin{equation}
\label{eq:4.21}
c>N^{\delta n} \log N^{ n}
\end{equation}
for some very large $n$, then the probability that  $G(N^{n}, \overline{c}/N^{n})$ is connected is close to $1$, by Theorem 2.8.1 in \cite{D}. It is possible to choose $c$ large enough so that (\ref{eq:4.21}) holds because the only restriction on $c$ is $c<N^{(1+\delta)n}$. This implies that the probability that all the points in  the $n$-ball are  connected is close to $1$. Then, taking $n=k_{n_0}$, (\ref{eq:4.10}) is true.

Finally, the uniqueness follows from   Theorem 2 in  \cite{KMT} (see Remark \ref{KMTr}).

 %we show that the percolation cluster constructed in this way is unique.  Suppose that there are two chains of nested balls $(B_{k_n})$ and $(B'_{k_n})$ %(containing two different given fixed points respectively) with their respective clusters $X_{k_n}$ and $X'_{k_n}$, each cluster being connected with the others %within its chain of nested balls for almost all $n$, forming an infinite cluster. As $n$ grows, $B_{k_n}$ and $B'_{k_n}$ will have a point in common for some %$n$ (since they grow to the left and to the right, think of the tree representation of $\Omega_N$).  Once $B_{k_n}$ and $B'_{k_n}$ have a point in common, by %ultrametricity they are the same ball, and since $|X_{k_n}|\geq N^{\gamma k_n}$ and $|X'_{k_n}|\geq N^{\gamma k_n}$ and $\gamma$ can be taken so that $N^{\gamma %k_n}>\frac{1}{2} N^{k_n}$ (see (4.4)), then $X_{k_n}$ and $X'_{k_n}$ cannot be disjoint, so they must coincide since they are maximal connected sets. From then %on there is only one infinite cluster. So the uniqueness of an infinite cluster  satisfying $|X_{k_n}|\geq N^{\gamma k_n}$ for $\gamma$ sufficiently close to % $1$ follows.
 \hfill $\Box$

%\textbf{We can also show that there cannot be more than one infinite cluster disjoint with %positive density: - is this correct?}

\begin{remark}\label{KMTr} \rm{
Uniqueness of the infinite cluster is proved in \cite{KMT}  (Theorem 2). They prove that $(\Omega_N,d)$ can be embedded into $\mathbb{Z}$ such that any ball of radius $r$ will be represented by $N^r$ consecutive integers, and the collection of balls of radius $r$ partitions $\mathbb{Z}$ and  that the embedding is  stationary and ergodic.  The uniqueness then follows from Gandolfi et al. \cite{GKN}, Theorem 0, on the uniqueness of the infinite cluster for long range percolation on $\mathbb{Z}$
satisfying the positive finite energy condition. The proof uses only the properties that the connection probabilities between vertices $\x,\y$  are strictly positive
and depend only on $d(\x,\y)$, and therefore is also applicable to our case.
An intuitive argument for uniqueness of the percolation cluster (of positive density), using the argument of the theorem, is that two chains of nested balls will eventually intersect, and by ultrametricity from then on they coincide, so, if their largest clusters occupy a sizeable part of the balls, they will eventually be the same.
}
\end{remark}

 %Our approach  takes advantage of the ultrametricity directly but does not exclude the possibility of ``thin'' (zero density) infinite clusters. See Section 6, %open problem 2.}

\section{Proofs for the critical case $\delta =1$}
\label{S5}
\setcounter{equation}{0}
\subsection{Lemmas for the proof of Theorem 3.5.}
\label{sub:5.1}

\begin{lemma}
%\label{L5.6}
\label{L5.1}
Let $b=0$ and with $aN$ instead of $a$ in (\ref{eq:2.5}), and $K=1$ in (\ref{eq:2.3}).

\noindent
(a) Let $A_{n,j}$ denote the event that a $k_n$-ball $B_{k_n}$ (containing ${\bf 0}$) is connected to the $(k_{n+1} +j-1, k_{n+1}+j]$-annulus around $B_{k_n}$, $j=1,\ldots ,k_{n+2}-k_{n+1}$. Then
\begin{equation}
\label{eq:5.5x}
 P(A_{n,j})\sim
{\frac{aN(1-N^{-1})\log n}{N^{j} N^{\log n}}} \quad \hbox{  as  } n\to\infty.
\end{equation}

\noindent
(b) Let $\widehat{A}_{n,j}$ denote the event that there in no connection between the  $k_{n+1}$-ball $B_{k_{n+1}}$ and the $(k_{n+1} +j-1, k_{n+1}+j]$-annulus around   $B_{k_{n+1}}$, $j=1, \ldots, k_{n+2} -k_{n+1}$. Then
\begin{equation}
\label{eq:5.6}
P(\widehat{A}_{n,j})\sim n^{-aN(1-N^{-1})/N^{j}}\quad {\it as}\quad n\rightarrow \infty .
\end{equation}

\noindent
(c) Let $A_n$ denote the event that there is no connection between the interior of a $k_{n+1}$-ball $B_{k_{n+1}}$ and the
$(k_{n+1}, k_{n+2}]$-annulus around $B_{k_{n+1}}$. Then
\begin{equation}
\label{eq:5.7}
P(A_n) \sim n^{-a}\quad {\it as}\quad n\rightarrow \infty,
\end{equation}
and this implies that

(i) if $a<1$, then with probability $1$ there are infinitely many pairs  $((k_n, k_{n+1}],\, (k_{n+1},k_{n+2}])$ of successive annuli that are not connected,

(ii) if $a>1$, then with probability $1$ there are at most finitely many pairs  $((k_n, k_{n+1}],\, (k_{n+1},k_{n+2}])$ of successive annuli that are not connected.

\noindent
(d) Let $A_{n,j,\ell}$ denote the event that  a $k_n$-ball $B_{k_n}$ is connected to the $(k_{n+\ell}+j-1, k_{n+\ell}+j]$-annulus around $B_{k_n}$, $j=1, \ldots, k_{n+\ell +1}-k_{n+\ell}$. Then for $\ell \geq 1$,
\begin{equation}
\label{eq:5.8}
P(A_{n,j,\ell})\sim \frac{aN(1-N^{-1}) \log (n+\ell)}{N^{j} (n+\ell)^{\ell \log N}}\quad {\it as}\quad n\rightarrow \infty.
\end{equation}
(Note that (\ref{eq:5.6}) is a special case with $\ell=1$).

\noindent
(e) Let $\widetilde{A}_n$ denote the event that there are  connections from a $k_n$-ball $B_{k_n}$ to the complement of a $k_{n+1}$-ball $B_{k_{n+1}}$, with $B_{k_n}\subset B_{k_{n+1}}$. Then
\begin{equation}
\label{eq:5.9}
P(\widetilde{A}_n) \lesssim\; a  \sum^\infty_{\ell =1}
\frac{\log (n+\ell)}{(n+\ell)^{\ell \log N}} \quad {\it as}\quad n\rightarrow \infty.
\end{equation}
(f) Let $A$ denote the event that there are no connections from a $k_n$-ball $B_{k_n}$ to the complement of a $k_{n+1}$-ball $B_{k_{n+1}}$, with $B_{k_n} \subset B_{k_{n+1}}$, for all but finitely many $n$. Assume in addition that in the case $N=2$, $K$ in (2.3) is taken so that $K>\frac{1}{\log 2}$. Then for any $a>0$,
\begin{equation}
\label{eq:5.10}
P(A)=1.
\end{equation}
\end{lemma}
\noindent
{\bf Proof.}

\noindent
(a) By (\ref{eq:2.5}) and (\ref{eq:2.6}), $c_{k_{n+1}+j}\sim c_{k_n}$ as $n\to\infty$, hence
\begin{eqnarray*}
P(A_{n,j}) &= & 1-\biggl( 1-
\frac{c_{k_{n+1}+j} }{N^{2(k_{n+1}+j)}}\biggr)
^{N^{k_n}N^{k_{n+1}+j} (1-N^{-1})}\\
&\sim & 1-\exp \biggl\{ -
\frac{aN(1-N^{-1})\log (n+1) }
{N^{j} N^{k_{n+1}-k_n}}\biggr\}\\
&\sim & 1-\exp \biggl\{-
\frac{aN(1-N^{-1})\log (n+1) }
{N^{j}N^{\log n}}\biggr\},
\end{eqnarray*}
and then (\ref{eq:5.5x}) follows.

\noindent
(b) Since there are $N^{k_{n+1}-k_n -1}\sim N^{\log n}\;\; k_n$-balls in the interior of the $k_{n+1}$-ball, then from (\ref{eq:5.5x}),
\begin{eqnarray*}
P(\widehat{A}_{n,j}) & \sim & \biggl(1-\frac{aN(1-N^{-1}) \log n}{N^{j}N^{\log n}}\biggr)^{N^{\log n}}\\
&\sim & \exp \biggl\{ -\frac{aN(1-N^{-1})\log n}{N^{j}}\biggr\}\\
&=& n^{-aN(1-N^{-1}) /N^{j}}.
\end{eqnarray*}
(c) By (\ref{eq:5.6}) and independence,
\begin{eqnarray*}
P(A_n) &\sim & \prod^{k_{n+2}-k_{n+1}}_{j=1} n^{-aN(1-N^{-1})/N^{j}}\\
&= & n^{-aN (1-N^{-1}) \sum^{\log (n+1)}_{j=1} N^{-j}}\\
&=& n^{-a (1-N^{-\log n})}\\
&\sim & n^{-a}.
\end{eqnarray*}
Since $\sum_n n^{-a} <\infty$ if and only if $a>1$, then by independence and the (second) Borel-Cantelli lemma, for $a<1$, with probability $1$ there are infinitely many successive $(k_n, k_{n+1}]$-annuli that are not connected, and  by the (first) Borel-Cantelli lemma, for $a>1$, with probability $1$ there are at most finitely many successive $(k_n, k_{n+1}]$-annuli that are not connected.

\noindent
(d) By (\ref{eq:2.5}) and (\ref{eq:2.6}),
\begin{eqnarray*}
P(A_{n,j, \ell}) &\sim & 1-\biggr(1-
\frac{C+aN\log (n+\ell)}{N^{2(k_{n+\ell}+j)}}\biggr)^{N^{k_n}N^{k_{n+\ell}+j}(1-N^{-1})}\\
&\sim & 1-\exp \biggl\{ -\frac{aN(1-N^{-1})\log (n+\ell)}
{N^{j}N^{k_{n+\ell}-k_n}}\biggr\}\\
&\sim& \frac{aN(1-N^{-1})\log (n+\ell)}{N^{j} N^{\ell \log (n+\ell)}}\\
& =&
\frac{aN(1-N^{-1})\log (n+\ell)}{N^{j} (n+\ell )^{\ell \log N}}.
\end{eqnarray*}
(e) By (\ref{eq:5.8}),
\begin{eqnarray*}
P(\widetilde{A}_n) &\sim &\sum^{\infty}_{\ell =1} \sum^{k_{n+\ell +1}-k_{n+\ell}}_{j=1}
\frac{aN(1-N^{-1})\log (n+\ell)}{N^{j} (n+\ell )^{\ell \log N}}\\
&\sim & aN \sum^\infty_{\ell =1}(N^{-1}-N^{-\log (n+\ell)}) \frac{\log (n+\ell)}{(n+\ell )^{\ell \log N}}\\
& \lesssim & a \sum^{\infty}_{\ell =1}\frac{\log (n+\ell)}{(n+\ell )^{\ell \log N}}.
\end{eqnarray*}
(f) By (\ref{eq:5.9}),
\begin{eqnarray*}
\sum^\infty_{n=1} P (\widetilde{A}_n) &\lesssim & a \sum^\infty_{n=1} \sum^\infty_{\ell =1} \frac{\log (n+\ell)}{(n+\ell )^{\ell \log N}}\\
&\lesssim & aN \sum^\infty_{j=2} \frac{\log j}{j^{\log N}}<\infty ,
\end{eqnarray*}
and the result follows by Borel-Cantelli.

This proof has been done for $K=1$ and $\log N>1$, hence $N\geq 3$.  For $N=2$ we have $k_{n+\ell}-k_n\sim K\log(n+\ell)$ (see (\ref{eq:2.4})), which yields $K\log N$ instead of $\log N$ in the last step, so we take $K>\frac{1}{\log 2}$ for summability.
 \hfill $\Box$

\begin{lemma}
\label{L5.7} Let $K$ and $b$ be as in part (b)(2) of Theorem \ref{T3.2}, that is, $0<b<2K-\frac{1}{\log N}$.
     Let $A_{n,j}$ denote the event that the cluster $X_{k_n}$ in a $k_n$-ball $B_{k_n}$ is connected to the $(k_{n+j}, k_{n+j+1}]$-annulus around $B_{k_n}$, $j\geq 2$. Then there is a positive constant $M$ such that
$$
P(A_{n,j}) \lesssim \frac{M\log n}{n^{(Kj-b)\log N}},
$$
as $n\to \infty$, and
$$\sum_{n=1}^\infty \sum_{j=2}^\infty P(A_{n,j})<\infty.$$
Hence with probability $1$ there exists a (random) number $n_0$ such that for all $n\geq n_0$ the connections between the clusters $X_{k_n}$ restricted to the $(k_n,k_{n+1}]$-annuli do not skip over two successive  annuli, that is, there are no connections between the annulus $(k_{n-1},k_n]$ and the annuli $(k_{n+2},k_{n+3}]$, $(k_{n+3},k_{n+4}]$, etc.
\end{lemma}
\noindent
{\bf Proof.}
By (\ref{eq:2.5}) and (\ref{eq:2.6}),

\begin{eqnarray*}
P(A_{n,j})
&=& \sum_{\ell =k_{n+j}+1}^{k_{n+j+1}}\left[1-\left(1-\frac{c_\ell}{N^{2\ell}}\right)\right]^{|X_{k_n}|N^\ell(1-N^{-1})}\\
&\leq& (C+a\log (n+1+j)\cdot N^{b\log (n+1+j)})N^{k_n}\sum_{\ell =k_{n+j}+1}^{k_{n+j+1}}\frac{1}{N^\ell}\\
&\lesssim & M\frac{\log n\cdot N^{b\log n}}{N^{k_{n+j}-k_n}}\lesssim M\frac{\log n}{N^{(Kj-b)\log n}}.
\end{eqnarray*}

%1-\biggl(1-\frac{c_{k_{n+j+1}}}{N^{2k_{n+j+1}}}\biggr)^{|X_{k_{n}}|(N^{k_{n+j+1}}-N^{k_{n+j}})}\\
%&<&
%\frac{c_{k_{n+j+1}}N^{k_n} N^{k_{n+j+1}}(1-N^{-(k_{n+j+1}-k_{n+j})})}
%{N^{2k_{n+j+1}}}\\
%&=& \frac{C+a\log (n+j+1) N^{b\log (n+j+1)}}{N^{k_{n+j+1}-k_n}}\\
%&< & M \frac{\log (n+j+1)}{(n+j+1)^{(K(j+1)-b)\log N}}.

It is easy to show that the assumptions on $K$ and $b$ imply that
$$
\sum_{n=1}^\infty\sum^\infty_{j=2} \frac{\log n}{n^{(Kj-b) \log N}}<\infty,
$$
and then the result  follows by Borel-Cantelli.

\hfill $\Box$
\vglue .25cm
  The main tools for proving part (b) of the theorem  are a large deviation inequality for the binomial distribution and a connectivity result for an E-R random graph. Recall that a graph is said to be connected if it has only one connected component and no isolated vertices.
\medskip

We first recall the large deviation bound for the binomial distribution  \cite{JLR} (Corollary 2.4).
\begin{lemma}
\label{L5.2} Let $Y_n$ be Bin$(n, 1-p)$ and $c>1$. Then for $x\geq cn(1-p)$,
$$
P(Y_n \geq x)\leq e^{-h(c)x},
$$
where $h(c)=\log c-1+1/c>0$.
\end{lemma}
\begin{corollary}
\label{C5.3}
There exist $\kappa >0$ and $\varepsilon >0$ such that for $0<\sigma < \min(\frac{p}{1-p},\varepsilon)$,
\begin{equation}
\label{eq:5.2}
P\biggl(1-\frac{Y_n}{n} \leq p-\sigma (1-p)\biggr)\leq e^{-\kappa \sigma^2 (1-p) n}.
\end{equation}
\end{corollary}
\noindent
{\bf Proof.} It is easy to see that there exist $\kappa >0$ and $\varepsilon>0$ such that $h(c)\geq \kappa (c-1)^2$ for $1<c<1+\varepsilon$. Then the result follows from the lemma putting $c=1+\sigma$. \hfill $\Box$

\begin{definition}{\rm
\label{D5.4} For $0<\beta <1$, we say that a $k$-ball $B_k$ is} $\beta$-good {\rm if its
   cluster $X_k$  satisfies $|X_k|\geq \beta N^k$.}
\end{definition}

Note that if ${\bf 0} \in B_k$ and $B_k$ is $\beta$-{\it good}, then  the probability that
${\bf 0} \in X_k$  is $\geq \beta$.

\begin{lemma}
\label{L5.5'}
 Assume that for some $\beta>0$ and some $n_0$,
\begin{equation}
\label{assumption}
P(|X_{k_n}|\geq \beta N^{k_n}\quad \hbox{\it for all}\quad n\geq n_0)=p_*>0.
\end{equation}
Then percolation occurs.
\end{lemma}
\noindent
{\bf Proof.} By transitivity we may assume that $\mathbf{0}$ belongs to the $k_n$-ball whose largest cluster is $X_{k_n}$. Then
 by the assumption (\ref{assumption}) we have that
$$\liminf_{k\rightarrow \infty} P({\bf 0}\in X_{k_n})\geq \beta p_*>0,
$$
which implies percolation. \hfill $\Box$

\begin{lemma}
\label{L5.5} Let  $0<b  <2K$ in (\ref{eq:2.5}) with  $k_n$ as in (\ref{eq:2.3}), and let $0<\beta <1$. Let $X_{k_n}$ and $X'_{k_n}$ be the largest clusters in  two (disjoint) $\beta$-{\it good} $k_n$-balls in a $k_{n+1}$-ball. Then
\begin{equation}
\label{eq:5.3'}
P(X_{k_n}\;{\rm and}\;X'_{k_n}\;\hbox{\it are connected within the}\; k_{n+1}\hbox{\it -ball})   \gtrsim r_n(\beta) \hbox{ as }n\to\infty,
\end{equation}
 where
\begin{equation}
\label{eq:5.3}
r_n (\beta) = \frac{\beta^2 a\log n}{N^{(2K-b)\log n}}.
\end{equation}
\end{lemma}
\noindent
{\bf Proof.}  By (\ref{eq:2.6}) and (\ref{eq:2.5}),  $k_n<d(X_{k_n},X'_{k_n})\leq k_{n+1}$, so
\begin{eqnarray*}
%\kern4cm &&
&& P( X_{k_n}\;{\rm and}\;{X'}_{k_n}\;\hbox{\it are connected within the}\; k_{n+1}\hbox{\it -ball})\\
 & & \geq  1 -\biggl(1-\frac{c_{k_{n}}}{N^{2k_{n+1}}} \biggr) ^{|X_{k_n}||{X'}_{k_n}|}\\
& & \geq 1-\biggl(1-\frac{c_{k_{n}}}{N^{2k_{n+1}}}\biggr)^{\beta^2 N^{2k_n}}\\
&&> 1 -\exp \biggl\{-\frac{c_{k_{n}}\beta^2}{N^{2(k_{n+1}-k_n)}}\biggr\}\\
&&\sim  1-\exp \biggl\{ -
\frac{\beta^2 (C+a \log n\cdot N^{b\log n})}{N^{2K\log n}}\biggr\}\\
&&\sim r_n(\beta).
\end{eqnarray*}
\hfill$\Box$

%\vglue .25cm
\noindent
\subsection{Proof of Theorem \ref{T3.2}}

\noindent
%(a) The proof is similar to the one for Theorem 3.1(a). For $j\geq k$,
%\begin{eqnarray*}
%\lefteqn{P(B_k\;\hbox{\rm  is connected to the}\; (j,j+1]\hbox{\rm -annulus around it}) }\\
%&=& 1- \biggl(1-\frac{C(j+1)^\alpha}{N^{2(j+1)}}\biggr)^{N^k N^j (N-1)}\\
%&\leq & M \frac{N^k(j+1)^\alpha}{N^{j+1}},\quad M\;\hbox{\rm is a constant},
%\end{eqnarray*}
%and since
%$$
%\sum_{k\geq 1}  \sum_{j\geq k}\frac{N^{k}j^\alpha}{N^j}<\infty \quad
%\hbox{\rm for}\quad \alpha <-1,
%$$
%the result follows by  Borel-Cantelli.

\bigskip

\noindent
(a)  (i)  Lemma \ref{L5.1}(f)  guarantees that with probability 1 there exists $n_0$ such that if  $n>n_0$ there are no connections between the $(k_n,k_{n+1}]$-annulus and the complement of $B_{k_{n+2}}$.  Moreover, if $a<1$, by Lemma \ref{L5.1}(c)(i) there are infinitely many $n$ such that the $(k_n,k_{n+1}]$-annulus and the $(k_{n+1},k_{n+2}]$-annulus are not connected.   This    implies that with probability 1 there there exists some $n\geq n_0$ such that there are no connections from $B_{k_{n+1}}$ to the exterior.
\medskip

\noindent
(ii) The pre-percolation statement  follows from Lemma \ref{L5.1}(c)(ii).

\bigskip
%\section{Proof of Theorem \ref{T3.2}}
%\label{S6}
%\setcounter{equation}{0}

\bigskip

\noindent
(b)
(1)  We begin by indicating the main ideas of the proof. We consider a sequence of nested balls $B_{k_n}$ (containing ${\bf 0}$) and  their largest  clusters $X_{k_n}$. Recall that each $k_{n+1}$-ball is comprised of $N^{k_{n+1}-k_n}$ (disjoint) $k_n$-balls.

At each stage we will focus on the subset of the $k_n$-balls in a $k_{n+1}$-ball that are $\beta_n$-{\it good} (i.e., $|X_{k_n}|\geq \beta_n N^{k_n}$, see Definition 5.5), where $(\beta_n)_n$ is a sequence of numbers in $(0,1)$ to be determined below. By construction the events that different $k_n$-balls are $\beta_n$-good are independent, and by transitivity they all have the same probability
 \begin{equation}
\label{eq:6.1}
p^G_n (\beta_n) = P(|X_{k_n}|\geq \beta_n N^{k_n}).
\end{equation}
Let ${\cal N}_n$ denote the number of $\beta_n$-good $k_n$-balls and recall (\ref{eq:5.3'}) and (\ref{eq:5.3}).
%\begin{equation}
%\label{eq:6.2}
%q_n (\beta_n) =P (X_{k_n}\;{\rm and}\; X'_{k_n}\; \hbox{\rm are connected within %the}\;\;k_{n+1}\hbox{\rm -ball}).
%\end{equation}
%This probability is the same for any two different clusters $X_{k_n}$ and $X'_{k_n}$.
Now we consider the E-R random graph $G({\cal N}_n, r_n (\beta_n))$ whose vertices are the ${\cal N}_n\;\beta_n$-good $k_n$-balls in the $k_{n+1}$-ball with connection probability $r_n(\beta_n)$.

The key idea of the proof is to establish that with probability $1$ there is a (random) number $n_0$ such that for $n\geq n_0$,
\begin{equation}
\label{eq:6.3}
G({\cal N}_n, r_n (\beta_n))\quad\hbox{\rm is connected},
\end{equation}
which implies that
\begin{equation}
\label{eq:6.4}
|X_{k_{n+1}}|\geq {\cal N}_n \beta_n N^{k_n}.
\end{equation}
We denote by $E_n$ the event
\begin{equation}
\label{eq:6.5}
E_n =\{G({\cal N}_n, r_n (\beta_n ))\quad \hbox{\rm is connected}\}
\end{equation}
and
\begin{equation}
\label{eq:6.6}
p^A_n (\beta_n)=P(E_n).
\end{equation}
In order to prove (\ref{eq:6.3}) for all large $n$, by Borel-Cantelli it suffices to show that
\begin{equation}
\label{eq:6.7}
\sum_n (1-p^A_n (\beta _n))<\infty.
\end{equation}
We denote by $F_n$ the event
\begin{equation}
\label{eq:6.10}
F_n =\{{\cal N}_n \geq (1-\varepsilon_n) p^G_n (\beta_n )N^{k_{n+1}-k_n}\},
\end{equation}
where $\varepsilon_n \in (0,1)$ and $p^G_n (\beta_n)$ is given by (\ref{eq:6.1}), and
\begin{equation}
\label{eq:6.9}
p^B_n (\beta_n, \varepsilon_n)=P(F_n),
\end{equation}

The next key idea is to choose a sequence of numbers $\varepsilon_n \in (0,1)$ of the form $\varepsilon_n=n^{-(1+\theta)}$, $\theta$ to be chosen below, with $\sum_n \varepsilon_n <\infty$ such that
\begin{equation}
\label{eq:6.8}
\sum_n (1-p^B_n (\beta_n, \varepsilon_n))<\infty.
\end{equation}

If both events $E_n$ and $F_n$ occur, then by (\ref{eq:6.4})
$$
|X_{k_{n+1}}|\geq (1-\varepsilon_n) p^G_n (\beta_n) \beta_n N^{k_{n+1}},
$$
so the $k_{n+1}$-ball is $\beta_{n+1}$-good with
\begin{equation}
\label{eq:6.11}
\beta_{n+1} = (1-\varepsilon_n)p^G_n (\beta_n) \beta_n.
\end{equation}
Therefore, since  $E_n$ and $F_n$ are independent (because $E_n$ is defined in terms of distance $k_{n+1}$, and $F_n$ in terms of distance $k_n$), then

\begin{equation}
\label{eq:6.12}
p^G_{n+1} (\beta_{n+1}) \geq p^B_n (\beta_n, \varepsilon_n) p^A_n (\beta_n).
\end{equation}
We will show that for sufficiently large values of  $C$ and $a$ there exists a sequence
$\beta_n$ such that
\begin{equation}
\liminf_n \beta_n>0.
\end{equation}
and
\begin{equation}
\label{eq:6.13}
\sum_n (1-p^G_n (\beta_n))<\infty,
\end{equation}
in order to obtain the results (\ref{eq:3.1}) and (\ref{eq:3.2}).

Since the quantities involved in the scheme described above are interdependent, we need  to overcome the interaction among them. We proceed as follows:
\begin{itemize}
\item[$\bullet$] We  set $\varepsilon_n=n^{-(1+\theta)}$ for some
$0<\theta <\frac{K\log N}{2}-1$ (recall that $K>\frac{2}{\log N}$),  hence
\begin{equation}
\label{eq:5.26}
K\log N >2(1+\theta).
\end{equation}
\item[$\bullet$] In Steps 1 and 2 below we will obtain  estimates
\begin{equation}
\label{eq:6.14}
 \sum_{n\geq n_0} \biggl(1-p^B _n \biggl(\beta_n, \frac{1}{n^{1+\theta}}\biggr)\biggr)<
\infty\quad \hbox{ if  }p^G_n(\beta_n)\geq\frac{1}{2}\hbox{  for all  }n\geq n_0,
\end{equation}
and
\begin{equation}
\label{eq:6.15}
\sup_{\beta \geq \frac{1}{5}} \sum_{n\geq n_0} (1-p^A_n (\beta))<\infty.
\end{equation}

We will then show that we can choose  $n_0$, $C$ and $a$ such that $\beta_n\geq \frac{1}{5}$ and $p^G_n(\beta_n)\geq \frac{1}{2} $ for all $n\geq n_0$.

\end{itemize}

 To complete the proof we proceed step by step. We first verify (\ref{eq:6.14}) and (\ref{eq:6.15}).

\noindent
Step 1. Assume that  $p^G_n(\beta_n)\geq \frac{1}{2}$ for all $n\geq n_0$. We first focus on one $k_{n+1}$-ball. It contains $N^{k_{n+1}-k_n}\;k_n$-balls. The probability that a $k_n$-ball is $\beta_n$-good is given by $p^G_n (\beta_n)$ (see (\ref{eq:6.1})). The number ${\cal N}_n$ of  $\beta_n$-good $k_n$-balls in the $k_{n+1}$-ball is ${\rm Bin} (N^{k_{n+1}-k_n}, p^G_n (\beta_n))$. From (\ref{eq:6.10}), (\ref{eq:6.9}) and Corollary \ref{C5.3},

\begin{eqnarray*}
1-p^B_n (\beta_n, \varepsilon_n) &=& P\biggl(\frac{{\cal N}_n}{N^{k_{n+1}-k_n}} <(1-\varepsilon_n) p^G_n (\beta_n)\biggr)\\
&=& P\biggl(1-\frac{N^{k_{n+1}-k_n}-{\cal N}_n}{N^{k_{n+1}-k_n}} <
p^G_n (\beta_n) -\sigma_n (1-p^G_n (\beta_n))\biggr)\\
&\leq & \exp \{-\kappa \sigma^2_n (1-p^G_n (\beta_n))N^{k_{n+1}-k_n}\},
\end{eqnarray*}
where
$$
\sigma_n =\varepsilon_n \frac{p^G_n(\beta_n)}{1-p^G_n (\beta_n)},
$$
hence
$$
1-p^B_n (\beta_n, \varepsilon_n) \leq  \exp\biggl\{-\kappa \frac{(p^G_n (\beta_n))^2}{1-p^G_n (\beta_n)} \varepsilon^2_n N^{k_{n+1}-k_n}\biggr\}.
$$
Since we assumed that $p:=p^G_n (\beta_n) \geq \frac{1}{2}$, then
$p^2 \geq \frac{1}{2} (1-p)$, so
\begin{equation}
\label{eq:6.16}
1-p^B_n (\beta_n, \varepsilon_n)\leq \exp\biggl\{- \frac{\kappa}{2}\varepsilon^2_n N^{k_{n+1}-k_n}\biggr\}.
\end{equation}
Using $\varepsilon_n=n^{-(1+\theta)}$ from (\ref{eq:6.16}) and (\ref{eq:2.4}),
\begin{equation}
\label{eq:6.17}
1-p^B_n \biggl(\beta_n, \frac{1}{n^{1+\theta}}\biggr)\leq \exp \biggl\{-
\frac{\kappa}{2}\frac{1}{n^{2(1+\theta)}} N^{K\log n}\biggr\} =
\exp \biggl\{-\frac{\kappa}{2} n^{(K\log N-2(1+\theta))}\biggr\}.
\end{equation}

We can then conclude from (\ref{eq:6.17}) and (\ref{eq:5.26}) that
$$
\sum_n \biggl(1-p^B_n \biggl(\beta_n, \frac{1}{n^{1+\theta}}\biggr)\biggr)<\infty,
$$ if $ p^G_n(\beta_n)\geq \frac{1}{2}$ for all $n\geq n_0$,
so (\ref{eq:6.14}) will be verified.

\bigskip
\noindent
Step 2. We  prove the result for a fixed $b$ satisfying
\begin{equation}\label{eq:3.2x} \frac{2}{\log N}<K<b<2K-\frac{1}{\log N},
\end{equation}
and then recalling Remark 2.2 observe that a simple coupling argument shows that the result remains true for all larger values of $b$.

Assume that $\beta_n\geq \frac{1}{5}$ for $n\geq n_0$. This assumption will be verified in Step 3.
 Define another constant
 \begin{equation}
 K_1=2K-b,
 \end{equation}
so that
\begin{equation}\label{old5.33}
\frac{1}{\log N}<K_1<K.
\end{equation}

%XXXXX - DELETE THIS
% Recalling (\ref{eq:6.5}), (\ref{eq:6.6}), by Lemma 5.8  we have, for large $n$,
%\begin{eqnarray*}
%p^C_n (\beta_n) && =  P({G}(
%\mathcal{N}_n,r_n(\beta_n) ) \rm{\;\;  is\; connected})\\
%&&\geq P(\mathcal{N}_n\geq N^{K_1\log n})
%P(G(N^{K_1\log n},r_n(\beta_n))\rm{\;\;is\; connected}).
%\end{eqnarray*}
%\begin{eqnarray}\label{eq:5.32}
%1-p^C_n (\beta_n) && \leq P({G}(
%N^{K_1\log n},r_n(\beta_n) ) \rm{\;\;  is\; not\;connected})\\
%&&\quad + 1- P(\mathcal{N}_n\geq N^{K_1\log n})\nonumber
%\end{eqnarray}
%because $N^{K_1\log n}<(1-\varepsilon_n)p^G_n(\beta_n)N^{K(k_{n+1}-k_n)}$ for large $n$, by Step 1.
%Since $b>K$, $K_1<K$,  there is a $c\in (0,1)$ such that
%\begin{eqnarray}\label{eq:5.33x}
%\sum_n (1-P(\mathcal{N}_n\geq N^{K_1\log n}))= \sum (1-p^B_n(\beta_n, c) <\infty.
%\end{eqnarray}
%XXXXXXXXXXXXXXXXXX

Recalling (\ref{eq:5.3}), (\ref{eq:6.5}) and conditioning on the event
$\{\mathcal{N}_n\geq N^{K_1\log n}\}$,  we have
$$ p^A_n(\beta_n)\geq P\left(G\left(\mathcal{N}_n,r_n(\beta_n)\right)\hbox{ is connected}|\;\mathcal{N}_n\geq N^{K_1\log n}\right)P(\mathcal{N}_n\geq N^{K_1\log n}),$$
and therefore
\begin{eqnarray}\label{eq:5.32}
&&1-p^A_n(\beta_n)\leq P\left(G\left(\mathcal{N}_n,r_n(\beta_n)\right)\hbox{ is not connected}|\;\mathcal{N}_n\geq N^{K_1\log n}\right)\\
&&\qquad\qquad\qquad\qquad  +1- P(\mathcal{N}_n\geq N^{K_1\log n}).\nonumber
\end{eqnarray}
Since $p^G_n(\beta_n)\geq\frac {1}{2}$ for large $n$, $\varepsilon_n\to 0$ and $K_1<K$ by (\ref{old5.33}), then, for large $n$,
$$
N^{K_1\log n}<(1-\varepsilon_n)p^G_n(\beta_n)N^{K\log n},$$
so, by Step 1,
\begin{equation}\label{eq:5.33x}
\sum_n(1-(P(\mathcal{N}_n\geq N^{K_1\log n}))<\infty.
\end{equation}

Using $N^{K_1\log n}\leq\mathcal{N}_n\leq N^{K\log N}$, assuming $\beta_n\geq\frac{1}{5}$ and taking $a>25K\log N$, and applying the inequality in the Appendix we obtain
\begin{eqnarray*}
&& P\left(G\left(\mathcal{N}_n,r_n(\beta_n)\right)\quad\hbox{\rm is not connected}|\;\mathcal{N}_n\geq N^{K_1\log n}\right) \\
&&\leq M[(\log n)^{13}n^{K_1\log N\cdot(1-a/25K\log N)}+n^{-K_1\log N}+e^{-L(\log n)^3 n^{2 K_1\log N}}],
\end{eqnarray*}
where $M$ and $L$ are positive constants.
Since $K_1\log N>1$, the sum of the second terms converges, and the sum of the third term also converges.  The sum of the first terms converges if
$$ a>25\left(K\log N+\frac{K}{K_1}\right)=: a_*.$$

Hence for $a>a_*$, together with (\ref{eq:5.32}), (\ref{eq:5.33x}) we have
$$ \sum_n(1-p^A_n(\beta_n))<\infty$$
uniformly for $\beta_n\geq \frac{1}{5}$.
which implies (\ref{eq:6.15}).

\noindent
Step 3. We must show that the assumptions  on $\beta_n$ and $p^G_n(\beta_n)$  used in steps 1 and 2 are self-consistent, that is, we can choose $n_0$, $\beta_{n_0}$, and $p^G_{n_0}(\beta_{n_0})$  such that
\begin{eqnarray}
&&\label{eq:5.33} p^G_{n}(\beta_{n})\geq \frac{1}{2} \hbox{   for all  }n\geq n_0,\\
&&\label{eq:5.34} \beta_n\geq\frac{1}{5} \hbox{   for all  }n\geq n_0.
\end{eqnarray}
\medskip

We proceed as follows.  Given $\theta$ satisfying (\ref{eq:5.26}) and recalling  (\ref{eq:6.14}), (\ref{eq:6.15}) we can choose $n_0$ such that the following  products satisfy:
\begin{eqnarray}
  &&\label{eq:5.35}\prod_{n\geq n_0}\biggl(1-\frac{1}{n^{1+\theta}}\biggr)\geq \biggl(\frac{4}{5}\biggr)^{1/3},\\
  &&\text{and for any  }k\in\mathbb{N}\nonumber\\
&&\label{eq:5.37}\prod_{n=n_0}^{n_0+k}p^B_n\biggl(\beta_n,\frac{1}{n^{1+\theta}}\biggr) \geq \biggl(\frac{4}{5}\biggr)^{1/3}\\&&\qquad\qquad \hbox{ provided that }p^G_n(\beta_n)\geq \frac{1}{2}\hbox{   for }n= n_0,\dots,n_0+k\nonumber\\
 &&\label{eq:5.36}\prod_{n= n_0}^{n_0+k}p^A_n(\beta_n) \geq \biggl(\frac{4}{5}\biggr)^{1/3}\\&&\qquad\qquad \hbox{ provided that }\beta_n\geq \frac{1}{5}\hbox{   for  }n= n_0\dots,n_0+k.\nonumber
\end{eqnarray}

Now choose
\begin{equation}\label{eq:3.n0} \beta_{n_0}\geq \frac{1}{2}\end{equation}
and
choose $C$ in (\ref{eq:2.5}) sufficiently large so that (see (\ref{eq:6.1}))
\begin{equation}
\label{eq:5.30}
p^G_{n_0} (\beta_{n_0})\geq \frac{2}{3}.
\end{equation}

Then by (\ref{eq:6.11}) and (\ref{eq:5.35})
\begin{eqnarray} &\beta_{n_0+1} =\left(1-\frac{1}{n_0^{1+\theta}}\right)p_{n_0}^G(\beta_{n_0})\beta_{n_0}
 \\ &\quad\;\geq  \biggl( \frac{4}{5}\biggr)^{1/3}\cdot \frac{2}{3}\cdot \frac{1}{2} >\frac{1}{5},\nonumber
\end{eqnarray}
and by (\ref{eq:6.12}), (\ref{eq:5.37}) and (\ref{eq:5.36}),  we  have
\begin{eqnarray}
&&\label{eq:5.41}p^G_{n_0+1}(\beta_{n_0+1})\geq p^B_{n_0}\biggl(\beta_{n_0},\frac{1}{n_0^{1+\theta}}\biggr)p^A_{n_0}(\beta_{n_0})\\
&&\geq \biggl(\frac{4}{5}\biggr)^{2/3}\geq \frac{1}{2}.\nonumber
\end{eqnarray}
Now assume that
\begin{equation}
\label{eq:5.42}
\prod_{\ell=1}^{k-1}p^G_{n_0+\ell}(\beta_{n_0+\ell})\geq \biggl(\frac{4}{5}\biggr)^{2/3},
\end{equation}
and $\beta_{n_0+\ell}\geq \frac{1}{5}$ for $\ell=1,\dots,k-1$.
\medskip

Then by (\ref{eq:6.12}),  (\ref{eq:5.37}) and (\ref{eq:5.36}) we  have
\begin{eqnarray}
&&\label{eq:5.41x}\prod_{\ell =1}^{k}p^G_{n_0+\ell}(\beta_{n_0+\ell})\\&&\geq \prod_{\ell =0}^{k-1} p^B_{n_0+\ell}\biggl(\beta_{n_0+\ell},\frac{1}{(n_0+\ell)^{1+\delta}}\biggr)p^A_{n_0+\ell}(\beta_{n_0+\ell})\nonumber\\
&&\geq \biggl(\frac{4}{5}\biggr)^{2/3}\geq \frac{1}{2}.\nonumber
\end{eqnarray}
Moreover by (\ref{eq:6.11}), (\ref{eq:5.35}), (\ref{eq:5.30}) and (\ref{eq:5.42}),
 \begin{eqnarray*}
&&\beta_{n_0+k}=\beta_{n_0}\prod_{j=1}^k\frac{\beta_{n_0+j}}{\beta_{n_0+j-1}}\\
&&=\beta_{n_0}\prod_{j=0}^{k-1}\biggl(1-\frac{1}{(n_0+j)^{1+\theta}}\biggr)\prod_{j=0}^{k-1}p^G_{n_0+j}(\beta_{n_0+j})\\
&&\geq \frac{1}{2}\cdot \biggl(\frac{4}{5}\biggr)^{1/3}\cdot  \frac{2}{3}\cdot   \biggl(\frac{4}{5}\biggr)^{2/3} >\frac{1}{5}.
\end{eqnarray*}

Therefore by induction we have
\begin{equation}\label{eq:5.43} \prod_{n\geq n_0+1}p^G_n(\beta_n)\geq\biggl(\frac{4}{5}\biggr)^{2/3}
\end{equation}
and
\begin{equation}\label{eq:5.44}\beta_n\geq \frac{1}{5} \text{ for all  }n\geq n_0.
\end{equation}

Step 4. For  $a>a_*$, $n_0$ and $C$ chosen above we then have
$$
\beta_n \geq \frac{1}{5}\quad{\rm and}\quad p^G_n (\beta_n)\geq \frac{1}{2}\quad\hbox{\rm for all}\quad n\geq n_0.
$$
This together with (\ref{eq:5.35}), (\ref{eq:5.36})  yields the estimate (\ref{eq:6.13}) and
\begin{equation}
\prod_{n\geq n_0} p^G_n(\beta_n)>0,
\end{equation}
which implies the assumption of Lemma \ref{L5.5'}.
Then
(\ref{eq:6.13}) implies (\ref{eq:3.2}) (for some $n_{00}\geq n_0$) and (\ref{eq:5.44}) implies (\ref{eq:3.1}). Percolation then follows from (\ref{eq:3.1}), (\ref{eq:3.2}) and Lemma \ref{L5.5'}.
The uniqueness of the infinite cluster again follows from \cite{KMT}, Theorem 2.

%To show uniqueness of the infinite cluster (with positive density)  we proceed analogously to the proof of Theorem 3.1(b).  Assume that there are two disjoint %infinite clusters with $\liminf_{n\to\infty} \frac{|X_{k_n}|}{N^{k_n}}\geq \beta >0$ and containing two different given points  respectively. By ultrametricity %they belong to the same $k_n$-balls for all large $n$.  Then the probability that there is no connection between the two infinite clusters in the annulus %$B_{k_{n}}\backslash B_{k_{n-1}}$ is

%\begin{eqnarray*}
%\lesssim\left( 1- \frac{a\log n\cdot N^{b\log n}}{N^{2n\log n}}\right)^{\beta^2N^{2n\log n}}\to 0\quad \hbox{as  }n\to\infty,
%\sim \left(1-\frac{\beta^2a\log n}{}\right)
%\end{eqnarray*}
%and therefore the probability that there is no connection between them is $0$. Thus we obtain a contradiction.

(2) The proof follows immediately from Lemma  \ref{L5.7}.

\medskip

\noindent
(c)   Consider the  case $0<b\leq \frac{2}{\log N}<K$ but modifying the model by replacing the actual connection probabilities at distances $k_n+1,\dots,k_{n+1}$ with the lower
bound $\frac{c_{k_{n}}}{N^{2k_{n+1}}}$. In this case the lower bound on the connection probabilities in Lemma \ref{L5.5} can be replaced by the upper bound

\[ P( X_{k_n}\;{\rm and}\;{X'}_{k_n}\;\hbox{\it are connected within the}\; k_{n+1}\hbox{\it -ball})\leq \widetilde r_n\quad \rm{as}\;\;n\to\infty,\]
where
\[ \widetilde r_n(\beta)=  \frac{a\log n}{N^{(K-2/\log N)\log n}}\frac{1}{N^{K\log n}}.\]
We can then consider the E-R graph $G(\mathcal{N}_n,\widetilde r_n(\beta))$.
Assuming that $\mathcal{N}_n$ is of order $ N^{K\log n}$, in this case by Erd\H{o}s-R\'enyi theory the resulting random graph has only of order $\log(N^{K\log n})$ good $k_n$-balls in in the largest connected component in the $k_{n+1}$-ball . This  would imply that the limiting density of the largest connected component  in the $k_n$-balls decreases to 0 as $n\to\infty$ so that percolation does not occur.

To make this precise consider $G(\mathcal{N}_n,\lambda/\mathcal{N}_n)$ with $0<\lambda <1$ and  $|\mathcal{C}_n|$  the size of the largest connected component. We have
\[ \alpha(\lambda) =\lambda-1-\log \lambda >0,\]
and there exists $n_0$ such that for $n\geq n_0$,   $\frac{a\log n}{N^{(K-2/\log N)\log n}}<\lambda$. Then, given $\varepsilon >0$, for $n\geq n_0$,
\[ P(|\mathcal{C}_n|\geq (1+\varepsilon)(\log \mathcal{N}_n)/\alpha(\lambda))\leq \mathcal{N}_n^{-(1+\varepsilon)}/\lambda,\]
(see \cite{D}, page 39).

We then have that the probability that there are more than $\frac{(1+\varepsilon)}{\alpha(\lambda)}\log \mathcal{N}_n$ good $k_n$-balls in the largest connected component in the $k_{n+1}$-ball for infinitely many $n$ is $0$.  Therefore there cannot be an infinite connected component (with positive density).

\hfill $\Box$

\subsection{Proof of Theorem \ref{MR-1}}
\noindent
(a)
Note that if $c_k=C_0+C_1\log k+ C_2k^\alpha$, with $\alpha >2$, then in Theorem 3.5 we can choose $b,K$ such that
\[ \frac{2}{\log N} <K<b< \frac{\alpha}{\log N}.\]
Then
\begin{eqnarray*}  &&c_{k_n}=C_0+ C_1\log\lfloor Kn\log n\rfloor+ C_2\lfloor K^\alpha n^\alpha(\log n)^\alpha\rfloor
\geq C+a\log n\cdot n^{b\log N}\end{eqnarray*}
for sufficiently large $C_0$ and $C_2$, where $C$ and $a$ are as in Theorem 3.5(b).
The proof follows then from (2.3), (2.5), (2.6) and the assumptions on $b$ in Theorem 3.5(b) and Remark 2.2.

\medskip

\noindent
(b) If $C_2=0$ and $C_1<N$, then $c_{k_n}= C_0+C_1 (\log K+ \log n+\log\log n)\leq \widetilde C_0+aN\log n$  for some $0<a<1$ and $\widetilde C_0> C_0$ (with $K=1$ if $N\geq 3$). The result then
follows from Theorem 3.5(a)(i).

\medskip

\noindent
(c) The existence of $C_*$ follows by the argument in \cite{KMT} (Theorem 1(b)) as follows. The expected number of edges from a given
vertex is (see (2.1))
\[ \sum_{k=1}^\infty(N-1)N^{k-1}p_{(k)} \leq \sum_{k=1}^\infty (N-1)N^{k-1}\frac{(C_0+ C_1\log k+C_2k^\alpha)}{N^{2k}}\]
which is less than $1$ for sufficiently small $C_0,C_1,C_2$. The result follows by coupling the largest connected cluster containing a given point with a subcritical branching process  (see e.g. \cite{JLR}, page 109).
\hfill $\Box$

\bigskip

\noindent
{\bf Appendix. Connectivity of a random graph}
\vglue .25cm
Consider the E-R random graph $G(n, \frac{a \log n}{n})$, $a>0$. Using a random walk approximation for cluster growth in a susceptible-infected-removed epidemic model, Durrett \cite{D} proves the known result that $P(G(n, \frac{a \log n}{n})$ is connected) $\rightarrow 1$ as $n\rightarrow \infty$ if $a>1$. Putting together the parts of the proof  one obtains (see p. 64) the lower bound for $a>1$,
\begin{eqnarray*}
\lefteqn{P\biggl(G\biggl(n, \frac{a\log n}{n}\biggr)\;\quad \hbox{\rm is connected}\biggr)}\\
&\geq& \biggl[\biggl(1-\frac{14 (a \log n)^{13}e^{(13a\log n)/n}}{n^a}\biggr)\biggl(1-\frac{1}{n^{2.1}}\biggr)\biggl(1-\frac{1}{n^2}\biggr)\biggr]^n \\ && \cdot (1-e^{-(\log n)^3/100})^{n(n-1)}.
\end{eqnarray*}
Then using the inequalities $1-x>e^{-2x},\; 0<x<0.7968$, and $1-e^{-x}<x,\; x>0$ it follows that for $a>1$,
\begin{eqnarray*}
\lefteqn{P\biggl(G\biggl(n, \frac{a\log n}{n}\biggr)\;\quad \hbox{\rm is not connected}\biggr)}\\
&\leq & M[(\log n)^{13}n^{1-a}+n^{-1}+\exp(-L(\log n)^{13}n^2)],
\end{eqnarray*}
where $M$ and $L$ are positive constants.

%The proof of our Lemma \ref{L5.1} begins with this inequality, with $N^{\log n}$ in place of $n$.

\bigskip
\noindent
{\bf Acknowledgment.} We thank an anonymous referee for comments that helped us to improve the paper.


\begin{thebibliography}{99}
\bibitem{AN} M. Aizenman and C.M. Newman (1986). Discontinuity of the percolation density in one-dimensional $1/|x-y|^2$ percolation models, {\it Comm. Math. Phys.} 107, 611-647.
\bibitem{AB}  R. Albert and A. L. Barab\'asi (1999). Emergence of scaling in
random networks, {\em Science} 286, 509-512.
%\bibitem{AS} S.M. Athreya and J.M. Swart, Survival of contact
%processes on the hierarchical group, {\it Probab. Theory Relat. Fields}, to
%appear.
\bibitem{BR} A.L. Barab\'asi and E. Ravaz (2003). Hierarchical organization in
complex networks, {\it Phys. Rev. E.}
67, 026112.
\bibitem{BB} I. Benjamini and N. Berger (2001).  The diameter of long-range
percolation clusters on finite cycles,
{\em Random Structures Algorithms} 19, No. 2 (2001), 102-111.
\bibitem{BLP} I. Benjamini, R. Lyons  and Y. Peres (1999). Group-invariant
percolation on graphs, {\em Geom. Funct. Anal.} 9, 29-66.
\bibitem{BS} I. Benjamini and O. Schramm (1996). Percolation beyond $Z^d$, many questions and a few answers. {\it Elect. Commun. Probab.} 1, 71-82.

\bibitem{Be} N. Berger (2002) Transience, recurrence and critical behavior for
long-range percolation, {\it Commun. Math. Phys.} 226 (2002), 531-558.
\bibitem{Bi} M. Biskup (2004). On the scaling of chemical distance in long-range
percolation models, {\em Ann. Probab.} 32, 2938-2977.
\bibitem{Bis} M. Biskup (2010). Graph diameter in long-range percolation, Random Structures and Algorithms, 39, 210-227.
%10
\bibitem{Bo} B. Bollob\'as (2001). Random Graphs, 2nd ed., Cambridge, 2001.

\bibitem{Bo1} B. Bollob\'as, S. Janson and O. Riordan (2007).  The phase
transition in inhomogeneous random graphs, {\it Random Structures and
Algorithms} 31, 3-122.
%13
\bibitem{Bo2} B. Bollob\'as and O. Riordan, Percolation (2006). Cambridge, 2006.

\bibitem{BEI} D.C. Brydges, S.E. Evans and J.Z. Imbrie (1992). Self-avoidung random walk on a hierarchical lattice in four dimensions, Ann. Probab. 20, 82-124.
\bibitem{CE} P. Collet and J.-P. Eckmann (1978). A renormalization group
analysis of the hierarchical model in statistical mechanics, {\it Lecture Notes
in Physics} 74, Springer-Verlag.
\bibitem{CGS} D. Coppersmith, D. Gamarnik and M. Sviridenko (2002). The diameter of a long-range percolation graph, {\it Random Structures and Algorithms} 21, 1-13.
%16
\bibitem{DG1} D.A. Dawson and L.G. Gorostiza (2007). Percolation in a
hierarchical
random graph, {\it Comm.  Stochastic Analysis} 1, 29-47.
\bibitem{DGW} D.A. Dawson, L.G. Gorostiza and A. Wakolbinger (2005). Degrees of transience and recurrence and hierarchical random walks, {\it Potential Analysis} 22, 305-350.
\bibitem{DGW1} D.A. Dawson, L.G. Gorostiza and Wakolbinger (2004). Hierarchical random walks, in ``Asymptotic Methods in Statistics'', 173-193, {\it Fields Institute Communications}, A.M.S.
%19
\bibitem{D}  R. Durrett  (2006). Random Graph Dynamics, Cambridge.

\bibitem{Dy}  F. J. Dyson (1969).  Existence of a phase-transition in a one-dimensional Ising ferromagnet,
Comm. Math. Phys. 12, 91.

%\bibitem{Gil} E.N. Gilbert (1959). Random graphs, Annals of Math. Statist. 30, 1141-1144.

\bibitem{GKN} A. Gandolfi, M.S. Keane and C.M. Newman (1992). Uniqueness of the infinite component in a random graph with applications to percolation and spin glasses, Probab. Th. Related Fields 92, 511-527.

\bibitem{G} G. Grimmett (1999). Percolation, Second Ed., Springer, 1999.



%\bibitem{GKM} G.R. Grimmett, M. Keane and J.M. Marstrand, On the
%connectedness of
%a random graph, {\it Math. Proc. Cambridge Phil Soc.} 96(1984), 151-166.

%\bibitem{HJ} B.M. Hambly and J. Jordan (2004). A random hierarchical lattice: the series parallel %graph and its properties, {\it Adv. Appl. Probab.} 36 (2004), 824-838.

\bibitem{JLR} S. Janson, T. {\L}uczak and A. Ruci\'nski (2000), Random Graphs,
Wiley-Interscience.


\bibitem{K} J. Kleinberg (2001). Small-world phenomena and the dynamics of
information, in
{\em Advances in Neural Information Processing Systems (NIPS)} 14, 431-438.
\bibitem{K2} J. Kleinberg (2006). Complex networks and decentralized search
algorithms,
{\em Proc. ICM, 2006,} Volume 3, 1019-1044, European Mathematical
Society, Zurich.

\bibitem{KMT} S. Koval, R. Meester and P. Trapman (2010). Long-range
percolation on the hierarchical lattice, arXiv PR1004, 1251..

\bibitem{NS} C.M. Newman and L.S. \ Schulman (1986). One dimensional
$1/|j-i|^{s}$
percolation models. The existence of a transition for $s\leq2$,
{\it Comm. Math. Physics} 104, 547-571.

%\bibitem{O} N. O'Connell, Some large deviation results for sparse
%random graphs, {\it Probab. Theory Relat. Fields.} 110 (1998), 277-285.

\bibitem{PSN} I. Pak and T. Smirnova-Nagnibeda (2000). On non-uniqueness of percolation
on nonamenable Cayley graphs, {\it C.R. Acad. Sci. Paris}, t. 330, Serie 1. 495-500.

\bibitem{RTV} R. Rammal, G. Toulouse and M.A. Virasoro (1986). Ultrametricity
for physicists, {\it Rev. Mod. Phys.} 58,  765-788.

%31
\bibitem{SF} S. Sawyer and J. Felsenstein (1983) Isolation by distance in a
hierarchically clustered population, J. Appl. Probab. 20, 1-10.



\bibitem{Sch} W.H. Schikhof (1984). Ultrametric Calculus. An Introduction to $p$-adic analysis, Cambridge.
%32
\bibitem{Sc} L.S. Schulman (1983). Long-range percolation in one dimension, {\it J.
Phys. A}
16, no. 17, L639-L641.

%\bibitem{ES} E.V. Slud, Distribution inequalities for the binomial law,
%{\it Ann.
% Probab.} 5 (1977), 404-412.
\bibitem{Si} Ya. Sinai, Theory of Phase Transitions: Rigorous Results, (1982). Pergamon Press.

\bibitem{SA} D. Stauffer and A. Aharony (1994).  Introduction to Percolation
Theory, Taylor \& Francis.

\bibitem{T} P. Trapman (2010). The growth of the infinite long-range percolation chuster, {\it Ann. Probab.} 38, 1583-1608.

\bibitem{TV} T. S. Turova and T. Vallier (2010). Merging percolation on $Z^d$ and classical
random graphs: phase transition, Random Structures and Algorithms 36, 185-217.

\end{thebibliography}
\end{document}